\documentclass[final]{siamltex}
\usepackage[colorinlistoftodos]{todonotes}
\usepackage{algorithm,algorithmic}
\usepackage[top=1in, bottom=1in, left=1.2in, right=1.2in]{geometry}
\usepackage{amsmath,amssymb}
\usepackage{graphicx,subcaption}
\usepackage{multirow}
\usepackage{float}
\usepackage{tikz,pgfpages}
\usetikzlibrary{shapes,arrows,trees,mindmap,decorations.pathreplacing}

\newtheorem{propos}{Proposition}

\newcommand{\define}{\stackrel{\text{def}}{=}}
\newcommand{\krylov}[3]{{\mathcal{K}}_{#3}\left\{#1,#2\right\}}
\newcommand{\Span}[1]{\text{Span}\{#1\}}
\newcommand{\normtwo}[1]{\| #1\|_{2}}
\newcommand{\normr}[1]{\| #1\|^\dagger}
\newcommand{\norm}[1]{\|#1 \|}
\newcommand{\bigO}{{\mathcal{O}}}

\DeclareMathOperator*{\argmin}{arg\,min}

\newcommand{\bSigma}{\boldsymbol\Sigma}

\newcommand{\bphid}{\boldsymbol{\phi}_d}
\newcommand{\bphii}{\boldsymbol{\phi}_i}
\newcommand{\bphi}{\boldsymbol{\phi}}

\newcommand{\balpha}{\boldsymbol\alpha}
\newcommand{\bbeta}{\boldsymbol\beta}
\newcommand{\bchi}{\boldsymbol\chi}
\newcommand{\beps}{\boldsymbol\varepsilon}
\newcommand{\bmu}{\boldsymbol\mu}

\newcommand{\bzero}{\textbf{0}}
\newcommand{\nsp}{N_\text{sp}}

\newcommand{\bb}{\mathbf{b}}
\newcommand{\bc}{\mathbf{c}}

\newcommand{\be}{\mathbf{e}}
\newcommand{\bbf}{\mathbf{f}}
\newcommand{\bp}{\mathbf{p}}
\newcommand{\br}{\mathbf{r}}

\newcommand{\bu}{\mathbf{u}}
\newcommand{\bv}{\mathbf{v}}

\newcommand{\bx}{\mathbf{x}}
\newcommand{\by}{\mathbf{y}}
\newcommand{\bz}{\mathbf{z}}

\newcommand{\bEta}{\boldsymbol\eta}

\newcommand{\bA}{\mathbf{A}}
\newcommand{\bB}{\mathbf{B}}
\newcommand{\bC}{\mathbf{C}}
\newcommand{\bD}{\mathbf{D}}
\newcommand{\bE}{\mathbf{E}}
\newcommand{\bF}{\mathbf{F}}

\newcommand{\bH}{\mathbf{H}}
\newcommand{\bI}{\mathbf{I}}
\newcommand{\bJ}{\mathbf{J}}
\newcommand{\bK}{\mathbf{K}}

\newcommand{\bM}{\mathbf{M}}

\newcommand{\bQ}{\mathbf{Q}}
\newcommand{\bR}{\mathbf{R}}
\newcommand{\bS}{\mathbf{S}}
\newcommand{\bT}{\mathbf{T}}
\newcommand{\bU}{\mathbf{U}}
\newcommand{\bV}{\mathbf{V}}
\newcommand{\bW}{\mathbf{W}}
\newcommand{\bY}{\mathbf{Y}}
\newcommand{\bZ}{\mathbf{Z}}

\usepackage{tikz,pgfpages}

\title{Fast algorithms for Hyperspectral Diffuse Optical Tomography}
\author{Arvind K. Saibaba \thanks{Department of Electrical and Computer Engineering, Tufts University $\{\text{arvind.saibaba,eric.miller}\}@\text{tufts.edu}$} \and  Misha Kilmer\thanks{Department of Mathematics, Tufts University $\text{misha.kilmer}@\text{tufts.edu}$ } \and Eric L. Miller\footnotemark[1] \and Sergio Fantini\thanks{Department of Biomedical Engineering, Tufts University $\text{sergio.fantini}@\text{tufts.edu}$}}
\begin{document}
\maketitle

\begin{abstract}

The image reconstruction of chromophore concentrations using Diffuse Optical Tomography (DOT) data can be described mathematically as an ill-posed inverse problem.  Recent work has shown that the use of hyperspectral DOT data, as opposed to data sets comprising of a single or, at most, a dozen wavelengths, has the potential for improving the quality of the reconstructions.   The use of hyperspectral diffuse optical data in the formulation and solution of the inverse problem poses a significant computational burden.  The forward operator is, in actuality, nonlinear. However, under certain assumptions, a linear approximation, called the Born approximation, provides a suitable surrogate for the forward operator, and we assume this to be true in the present work. Computation of the Born matrix requires the solution of thousands of large scale discrete PDEs and  the reconstruction problem, requires matrix-vector products with the (dense) Born matrix.    In this paper, we address both of these difficulties, thus making the Born approach a computational viable approach for hyperspectral DOT (hyDOT) reconstruction. In this paper, we assume that the images we wish to reconstruct are anomalies of unknown shape and constant value, described using a parametric level set approach, (PaLS)~\cite{aghasi2011parametric}  on a constant background.  Specifically, to address the issue of the PDE solves, we develop a novel recycling-based Krylov subspace approach that leverages certain system similarities across wavelengths.   To address expense of using the Born operator 
in the inversion, we present a fast algorithm for compressing the Born operator that locally compresses across wavelengths for a given source-detector set and then recursively combines the low-rank factors to provide a global low-rank approximation. This low-rank approximation can be used implicitly to speed up the recovery of the shape parameters and the chromophore concentrations. We provide a detailed analysis of the accuracy and computational costs of the resulting algorithms  and demonstrate the validity of our approach by detailed numerical experiments on a realistic geometry.  
     
\end{abstract}

\section{Introduction}
Diffuse optical tomography (DOT) is an imaging technique that uses near infrared light to image highly scattering media. A good review has been provided in~\cite{arridge1999optical} and an updated version is provided in~\cite{arridge2009optical}. The imaging modality has shown great promise as a low-cost alternative or complement to existing medical imaging technology particularly in brain imaging and breast cancer detection.  The region of interest is illuminated with near infrared light over a collection of wavelengths and the data are comprised of observations of the resulting scattered diffuse fields at a number of locations surrounding the medium.  Given these measurements as well as the partial differential equation governing the interaction of light and tissue (typically, the diffusion equation), we seek to recover space and time-varying maps (i.e. images) of concentrations of physiologically relevant chromophores such as oxygenated and deoxygenated hemoglobin (HbO$_2$ and HbR respectively), lipid, and water (H$_2$O) as well as properties governing the scattering of light within the medium.

The recovery of images of chromophore concentrations can be mathematically posed as a nonlinear inverse problem. However, due to the diffusive physics associated with this problem as well as limitations concerning the geometric distribution of sources and detectors, image recovery is an ill-posed inverse problem. New technology developed in our research group allows for the collection of hyperspectral data (over 100 bands in the near infrared portion of the spectrum).  Although we have demonstrated~\cite{larusson2013parametric,larusson2011hyperspectral,larusson2012parametric} that the availability of more information using multiple wavelengths increases the accuracy of the reconstruction, the use of hyperspectral data poses a significant computational burden in the context of image recovery. We are interested in developing computationally efficient methods for hyperspectral diffuse optical tomography (HyDOT) with specific application towards breast imaging in which the breast is placed in between two parallel plates.

\begin{figure}[!ht]
\centering
\includegraphics[scale=0.35]{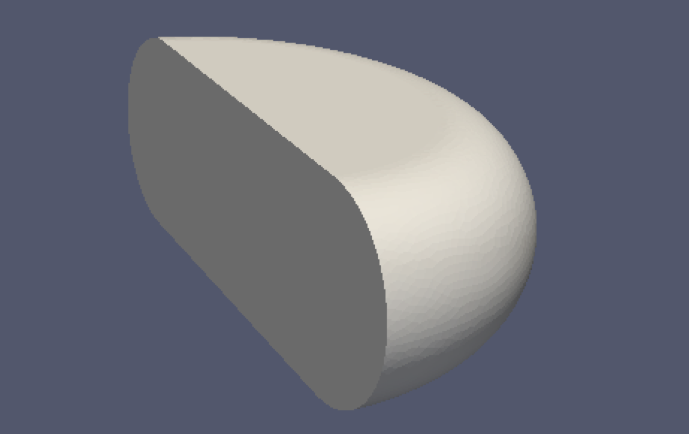}
\caption{The breast shaped phantom used as the imaging medium. At its widest, the phantom is $12$ cm long and $5$ cm thick.}
\label{fig:phantom}
\end{figure}

To motivate the need for fast algorithms for hyperspectral DOT, we outline here the expected costs in terms of storage and computation. To make ideas concrete, we consider an experimental setup for detecting tumors in breast tissue (see Figure~\ref{fig:phantom}). We use $N_s$ near-infrared sources to illuminate the medium of interest. The sources are constrained to lie on the top plane and detectors are constrained to be on a different plane so that for a given source, we have $N_{ds}$ detectors measuring photon fluence at $N_\lambda$ wavelengths. This results in $M = N_sN_{ds}N_\lambda$ measurements. We also assume that the domain has been discretized into a grid with $N$ vertices.

\begin{table}[!ht]
\centering
\begin{tabular}{|l|c|c|}\hline
Number of  & Symbol & Typical number \\ \hline
Sources & $N_s$ & $10-100$ \\ 
Detectors / source & $N_{ds}$ & $3-10$ \\ 
Wavelengths & $N_\lambda$ & $10-200$ \\ 
Grid size &  $N$ & $32^3-100^3$ \\  \hline
\end{tabular}
\caption{Typical range of parameters for the hyperspectral DOT problem}
\label{tab:range}
\end{table}
In Table~\ref{tab:range}, we have listed the range of various parameters that one might encounter in practice. Since the most accurate forward model in terms of the unknown voxel values is nonlinear, standard numerical/optimization approaches to solve the inverse problem repeatedly linearize the forward problem about a current estimate~\cite{arridge1999optical,de2011regularized}. Each optimization step then requires the solution of the forward and the adjoint PDE for each source-detector set~\cite{arridge1999optical}. In all, we need to solve $N_s(N_{ds} + 1)N_\lambda$ systems of equations \textit{at each optimization step}, which amounts to about $2\times 10^5$ systems of equations for the range of parameters described in Table~\ref{tab:range}. Thus, for finely discretized fields, even with a solver of optimal complexity $\bigO(N)$, use of the nonlinear forward model poses a significant computational challenge because the resulting cost is $\bigO(N_s(N_{ds} + 1)N_\lambda N$) flops. 

To mitigate the computational burden, in the present work we restrict ourselves to the case where the absorption can be represented as a small perturbation about the background medium.  Then, we can use the Born approximation to linearize the forward model.  Although this approximation introduces a modeling error and has known limitations~\cite{boas1997fundamental}, studies with experimental data have validated the utility of this model for hyperspectral DOT~\cite{larusson2013parametric,larusson2011hyperspectral,larusson2012parametric}. However, since the domain of interest has a complex shape, the Green's function required for the construction of the linearized operator is not known in closed form.  Therefore, we need to use a PDE formulation to compute the incident and adjoint fields.   An advantage of the PDE formulation is that we can handle known, non-constant background which may be obtained by imaging using a high resolution modality (such as magnetic resonance imaging, MRI) along with image segmentation~\cite{boverman2005quantitative}. Although the breast is a highly heterogenous medium, based on previous work we model it as a piecewise constant medium with homogenous background. For example, the authors in~\cite{schweiger1999optical} use prior anatomical knowledge to derive a piecewise constant medium.  Although the Born approximation reduces the computational burden associated with the imaging problem tremendously, computing the Born operator $\bH$ (which is a linear mapping between the perturbation of interest and the measurements) still 	requires the solution of many large-scale discretized PDEs for the incident and adjoint field; in sum, about $\bigO(N_s(N_{ds} + 1)N_\lambda)$ large scale, sparse linear systems corresponding to the discrete diffusion equation.
 
Furthermore, the cost to store $\bH$ and form matrix-vector products (matvecs) with $\bH$ during the optimization that are needed to invert for the desired parameters scales as $\bigO(MN)$, where $M$ is the number of measurements and $N$ is the number of unknowns.  Therefore, we develop a recursive algorithm to generate a low-rank approximation to $\bH$ and use this approximation in the optimization.  If the numerical rank of the low-rank representation is $R$ then the cost of storing factors and forming matvecs are $\bigO(R(M + N))$.
Of course, the optimal rank-R approximation could in theory be computed from the SVD of $\bH$ if we had $\bH$. This is inefficient on two fronts: a) it requires the full computation and storage of as well as multiple accesses to the very large, dense matrix $\bH$ b) the cost of an SVD on $\bH$, if we had it, is still prohibitively expensive at $\bigO \left( \min \{M,N\}^2, \max\{M,N\} \right)$ flops.   On the other hand, our algorithm, because it operates in a spatially recursive manner requiring, at the finest scale, only local information, and using rank revealing factorizations to aggregate information globally, does not require that $\bH$ be fully stored.

We have three main goals in this paper.  The first is to design an algorithm to overcome the challenge of computing the solution of $\bigO(10^5)$ number of large-scale parametric linear systems needed to obtain certain information necessary to compute our approximate Born matrix.   The second is to design an efficient algorithm to compute and store a meaningful low-rank representation of the measurement operator $\bH$.   The third
is to use this low-rank representation for recovering the  parameters that define our chromophore image.  We now summarize the key components of our fast algorithm and contributions in this paper:
\begin{itemize}

\item For each source-detector set, we need to compute the incident field and the adjoint field for hundreds of wavelengths.   We derive a novel Krylov recycling subspace approach to solve the corresponding systems of large-scale parametric linear systems which takes advantage of similarities in the systems across wavelengths.  This is described in detail in Section~\ref{sec:krylov}.

\item  For the problem at hand, the sources, detectors and the perturbation are well-separated from each other. Under these circumstances there is strong numerical evidence that the Born operator $\bH$ is low-rank, a feature that has been exploited to achieve computational savings~\cite{chaillat2012faims}.  In this work, we give a new approach for computing a low rank approximation, and use this in our numerical results on the parametric inverse problem.  
The storage of the matrix $\bH$ in its entirety is avoided; instead, the rows of $\bH$ corresponding to a single source-detectors set are constructed, compressed and then recursively compressed across multiple source-detectors pairs to obtain a low-rank factorization  $\bH \approx \hat{\bH} = \bU \bV^T$. We provide a detailed outline of the algorithm and analysis of the computational costs  in Section~\ref{sec:compress}. Because we fully compute the entries of $\bH$,  we need the fast Krylov solvers because we need to solve $\bigO(10^5)$ systems. 

\item   Based on recent success in the context of DOT, we employ the parametric level set (PaLS) approach to represent the chromophore image we want to recover.  The PaLS approach, developed for inverse problems in~\cite{aghasi2011parametric} and subsequently applied to diffuse optical tomography~\cite{larusson2012parametric,larusson2013parametric}, has the advantage of explicitly describing the geometry of the anomaly.

 As a result, the recovery of the chromophore image is obtained by solving a nonlinear least-squares solver on a problem in which we have replaced the Born operator by the approximation described above.  We show how to accelerate the reconstruction of the chromophore concentration and shape parameters by using the compressed measurement operator $\hat{\bH}$ is used in place of the full measurement operator. Error bounds are provided based on the error in the low-rank representation and are useful for the optimization routine. This is outlined in Section~\ref{sec:recon}. 
\end{itemize}

\textbf{Related work}: We briefly review other approaches to solve the parametric system of equations that represent the propagation of light in diffuse medium. One approach to deal with parametric system of equations is the use of spectral methods by expanding the matrix coefficients in a series of orthogonal polynomials (such as Chebyshev, Legendre, etc.) and solving for the coefficients of the orthogonal polynomials using a Collocation method or a residual minimizing Galerkin projection method (see~\cite{constantine2010spectral} and references therein). The low-rank property of the solutions arising from system of parametric coefficients with coefficients smoothly dependent on the parameters was demonstrated by Kressner and Tobler~\cite{kressner2011low} and they developed global Krylov subspace algorithms that exploited the low-rank nature to minimize computational and storage costs. In both approaches, one has to solve a coupled system of equations. This can be computationally expensive and, in order to ensure rapid convergence, a preconditioner that is effective across all the shifts is necessary but difficult to choose in practice. Other approaches to deal with the expensive cost of solving parametric linear systems is to use parametric model reduction which is reviewed in~\cite{benner2013survey}. 

Other works have also considered the compression of the measurement operator $\bH$ that maps the perturbation to the measurements. In~\cite{markel2003inverse}, the authors develop analytical formulas for inversion based on Fourier analysis when the sources and detectors are distributed uniformly on the boundary of a regular geometry such as a plane, cylinder, or sphere. In our previous work, we have exploited the structure of the Green's function in regular geometries to decompose the Born operator into a number of sparse easily computed matrices~\cite{hyde2007analysis}. The approach of compressing the operator $\bH$ is similar to that derived in~\cite{chaillat2012faims}. Here, the authors consider compression across multiple sources and detectors for a given frequency using randomized SVD and then recursively compressing the low-rank factors across multiple frequencies.  However, we cannot adopt their framework directly for the following reasons. Firstly, in our system, the detectors are not shared across all of the sources. Consequently, the pre-processing step that compresses the incoming field and the data, as described in~\cite{chaillat2012faims}, cannot be used directly.  Secondly, the authors in~\cite{chaillat2012faims} advocate compression of the measurement operator across multiple sources and detectors for a given frequency and then recursively combining the low-rank factorizations across different wavelengths. Memory limitations prevent computations of the entire measurement operator in its entirety. Because of the way our computations are organized, we choose to first compress across multiple wavelengths and detectors for a given source and then combine the low rank factorizations across multiple sources.

\section{Forward Problem}
In this section, we give the specifics of the forward problem and discretization used for the associated PDEs.

\subsection{Governing equations}
The radiative transport physics associated with the propagation of light through a medium can be approximated by the diffusion model of the form in the domain $\Omega$~\cite{arridge1999optical}
\begin{align}\label{eqn:diffusion} 
-\nabla \cdot D^\text{tot}(\br,\lambda) \nabla \phi(\br,\lambda) + \nu \mu_a^\text{tot}(\br,\lambda) \phi(\br,\lambda)  = & \quad   S(\br,\lambda)  & \quad \br \in& \Omega \\ \label{eqn:dirichlet}  
\phi(\br, \lambda) = &  \quad 0 & \quad \br \in& \partial\Omega_D\\ \label{eqn:robin} 
\phi(\br,\lambda) + 2AD^\text{tot}(\br,\lambda)\frac{\partial \phi(\br,\lambda)}{\partial n}  = &  \quad 0 & \quad \br \in& \partial\Omega_R 
\end{align} 
where $D^\text{tot}(\br,\lambda)$ is the diffusion coefficient and is related to the reduced scattering coefficient $\mu_s'(\lambda)$ as $D(\br,\lambda) = \nu/3\mu_s' (\br,\lambda) $.  We also denote by $\partial\Omega_D$ the portion of the boundary over which zero Dirichlet boundary conditions are imposed (curved boundaries and chest wall) and $\partial \Omega_R$ corresponds to the boundary portion over which Robin boundary conditions are imposed corresponding to a refractive index mismatch (top and bottom flat regions). The coefficient $A$ is a function of the refractive index of the medium. We denote by $\phi(\br,\lambda)$  the photon fluence at a position $\br$ due to a source of wavelength $\lambda$ injected into the medium,  and $\nu$ is the electromagnetic propagation velocity within the medium. Further, $\mu_a^\text{tot}(\br,\lambda)$ is the absorption coefficient. The quantity $S(\br,\lambda)$ is the photon source with units of optical energy per unit time per unit volume and typically written in terms of a delta function; that is, $S(\br,\lambda) = S_0(\lambda) \delta(\br-\br_s)$, with $S_0(\lambda)$ the source power at wavelength $\lambda$.

We decompose the absorption, $\mu_a^\text{tot}(\br,\lambda)$, into a constant background absorption $\mu_a(\lambda)$ and a spatially varying perturbation $\Delta \mu_a(\br,\lambda)$. The total fluence, $\phi$, is decomposed into an incident field $\phi_i$ and a scattered field $\phi_s$. Likewise, we can expand the diffusion $D^\text{tot}(\br,\lambda)$ into the sum of a homogenous background term $D(\lambda)$ and a perturbation $\Delta D(\br,\lambda)$. However, the spatial dependence of diffusion is minimal in healthy breasts~\cite{shah2004spatial} and spatial contrast in breast tumors is either non-existent or small~\cite{grosenick2004concentration,grosenick2005timea}. To simplify matters we assume that $D^{tot}(\br,\lambda)$ is independent of $\br$ and  $\Delta D(\br,\lambda)  =0$. Therefore $D^\text{tot}(\br,\lambda)$ can be represented entirely as $D(\lambda)$ and  we can then divide throughout by $D(\lambda)$ (see for example~\cite{larusson2011hyperspectral}). 
The equation for the incident field $\phi_i$ and the scattered field $\phi_s$ can be obtained by substituting $\phi(\br,\lambda) = \phi_i(\br,\lambda) + \phi_s(\br,\lambda)$ and collecting the appropriate terms, and is therefore,
\begin{align} 
\label{eqn:incident}
-\nabla^2 \phi_i(\br,\lambda) +  \frac{\nu \mu_a(\lambda)}{D(\lambda)} \phi_i(\br,\lambda) \quad  = & \quad  \frac{S_0(\lambda)}{D(\lambda)}\delta(\br-\br_s)  & \quad  \br \in \Omega \\ 
\label{eqn:born}
-\nabla^2 \phi_s(\br,\lambda) +  \frac{\nu \mu_a(\lambda)}{D(\lambda)} \phi_s(\br,\lambda) \quad = & \quad  - \frac{\nu\Delta\mu_a(r,\lambda)}{D(\lambda) }\left( \phi_i(\br,\lambda) + \phi_s(\br,\lambda) \right)   &  \quad \br \in \Omega  
\end{align}
Under the Born approximation, the scattered field is assumed to be much smaller than the incident field, i.e. $\phi_s  \ll \phi_i$ and therefore, the total fluence $\phi(\br,\lambda) = \phi_i(\br,\lambda) + \phi_s(\br,\lambda)$ in the right hand side of equation~\eqref{eqn:born} 
can be replaced by $\phi_i(\br,\lambda)$. 
As a result of this approximation, there is a linear relation between the scattered fluence rate $\phi_s(\br,\lambda)$ and the perturbation of absorption $\Delta\mu_a(\br,\lambda)$.  

 It should also be noted that both the scattered field $\phi_s(\br\,\lambda)$ and $\phi_i(\br,\lambda)$  satisfy the same boundary conditions in the equations~\eqref{eqn:dirichlet}-\eqref{eqn:robin}. Furthermore, if additional information such as spatial variability is known about the background properties of diffusion and absorption (currently assumed to be homogenous) they can be incorporated into this model~\cite{boverman2005quantitative}. The solution to the photon fluence $\phi_s$ computed at the measurement location $\br_d$ for a particular wavelength $\lambda$ can be written using the following integral equation 
\begin{equation}
\phi(\br_d, \lambda) =  \phi_i(\br_d, \lambda) + \phi_s(\br_d,\lambda) \quad \approx \quad \phi_i(\br_d,\lambda) - \int_\Omega \phi_d(\br,\lambda) \nu\Delta\mu_a(\br,\lambda) \phi_i(\br,\lambda)d\br \label{eqn:lippschwin}
\end{equation} 
where $\phi_d (\br,\lambda)$, which we call the adjoint field, can be derived using the reciprocity property of the Greens function and satisfies the system of equations along with the same boundary conditions in equations~\eqref{eqn:dirichlet}-\eqref{eqn:robin}
\begin{equation}
-\nabla^2 \phi_d(\br,\lambda) +  \frac{\nu \mu_a(\lambda)}{D(\lambda)} \phi_d(\br,\lambda) \quad  =  \quad  \frac{1}{D(\lambda)}\delta(\br-\br_d)   \qquad  \br \in \Omega
\end{equation}
and $\br_d$ corresponds to the detector location. 
To relate the scattered fluence to the concentrations of chromophores, the perturbation $\Delta\mu_a(\br,\lambda)$ is decomposed in terms of piecewise constant functions as 
\begin{equation} \Delta\mu_a(\br,\lambda) \define \sum_{l=1}^{N_{sp}} \varepsilon_l(\lambda)c_l\chi(\br)\qquad \chi(\br) = \left\{ \begin{array}{ll} 1 & \br \in \mathcal{D} \\ 0 & \br \in \Omega\backslash \mathcal{D}\end{array} \right. \label{eqn:perturbation}
\end{equation}
where $N_{sp}$ is the number of species, $\varepsilon_l$ is the extinction coefficient of species $l$ at wavelength $\lambda$, $c_l$ is the concentration of species $l$ and $\chi$ is an indicator function which depends on  $\mathcal{D}$, the domain of support for the perturbation we wish to image. For the purpose of this paper, we will consider that the chromophore concentrations are co-located. This choice was also considered in~\cite{larusson2012parametric}. Further details regarding the governing partial differential equations can be obtained from the following references~\cite{larusson2013parametric,larusson2011hyperspectral,larusson2012parametric}.

\begin{figure}[!ht]
\centering
\includegraphics[scale=0.29]{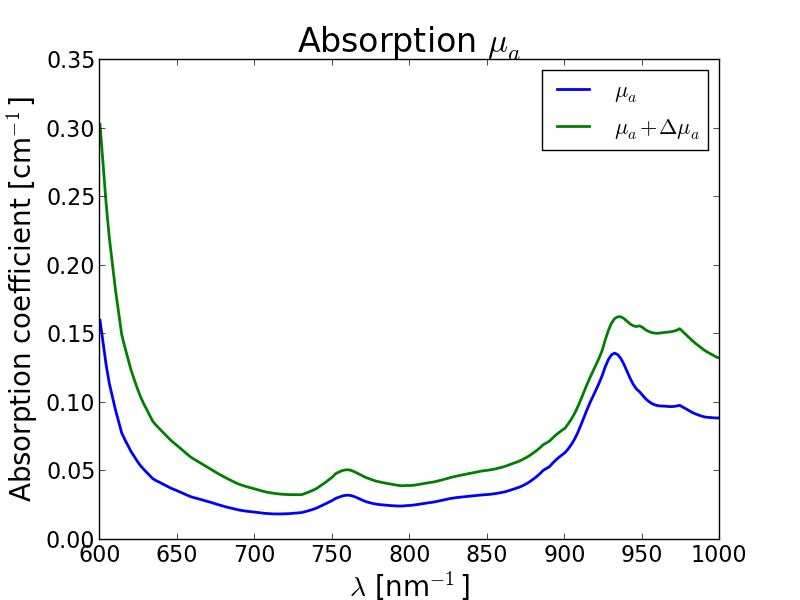}
\includegraphics[scale=0.29]{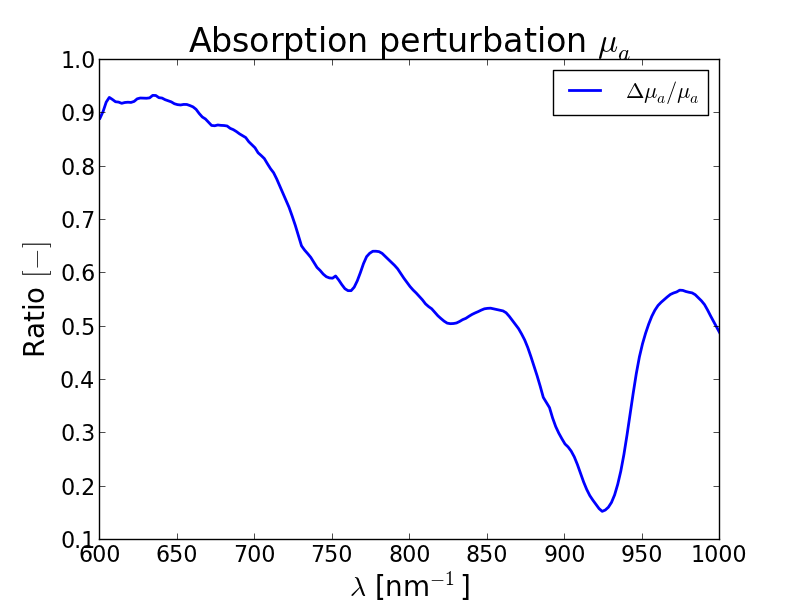}
\caption{The absorption coefficient $\mu_a$ as a function of $\lambda$ in the range $\lambda \in [600,1000]$ [nm]. The chromophore concentrations of the background and the perturbation used to generate this plot are provided in Table~\ref{tab:conc}.}
\label{fig:muaD}
\end{figure}
We also model the diffusion coefficient $D(\lambda)$ using Mie scattering theory~\cite{larusson2011hyperspectral} as 
\begin{equation}
D(\lambda) =\frac{\nu \Psi}{3} \left(\frac{\lambda}{\lambda_0}\right)^b
\label{eqn:difflambda}
\end{equation}
The reference wavelength $\lambda_0$ is chosen as $600$ nm and $\Psi$ has units of cm$^{-1}$. The scattering pre-factor $\Psi$ depends primarily on the number and size of scatterers, and a scattering exponent $b$ depends on the size of scatterers in the medium~\cite{grosenick2005timeb}.

\subsection{Discretization using finite elements}
To solve the systems of equations~\eqref{eqn:born} and~\eqref{eqn:incident} we use the standard linear Galerkin finite element approach. Expanding the solution field in an appropriately chosen finite dimensional basis $u_m(\br)$, i.e., \[ \phi_i(\br,\lambda) \approx \sum_{m=1}^N \hat{\phi}_{i,m}(\lambda) u_m(\br)\quad \text{and}\quad \phi_d(\br,\lambda) \approx \sum_{m=1}^N \hat{\phi}_{d,m}(\lambda) u_m(\br)\] Let us denote the discretized incident photon fluence field as $\bphii(\lambda) = [\hat{\phi}_{i,1}(\lambda),\dots,\hat{\phi}_{i,N}(\lambda)]^T$ and the scattered photon fluence $\bphi_d(\lambda) = [\hat{\phi}_{d,1}(\lambda),\dots,\hat{\phi}_{d,N}(\lambda)]^T$. The resulting system of equations can be summarized as 
\begin{equation} \label{eqn:discincident}
\left(\bK + \frac{\nu\mu_a}{D}(\lambda)\bM + \frac{1}{2AD(\lambda)} \bR\right) \bphi_i(\lambda)  =  \frac{1}{D(\lambda)}\bb_{i}
\end{equation} 
where the matrices $\bK$, $\bM$ and $\bR$  have entries given by 
\begin{align} \label{eqn:matrices}
\bK_{jk} \quad & =  \quad \int_\Omega \nabla u_k(\br) \cdot \nabla u_j(\br) d\br  \qquad \bR_{jk} = \int_{\partial\Omega_R}  u_j(\br)u_k(\br)d\br \\ \nonumber
\bM_{jk} \quad & =  \quad \int_\Omega  u_k(\br) u_j(\br) d\br 
\end{align}
for $j,k = 1,\dots,N$ and the vector $\bb_i$ has entries $\bb_{i,j} = \int_\Omega \delta (\br-\br_s) u_j(\br) d\br$. A similar equation can be derived for the adjoint field $\bphid$ with the same matrices and a different right hand side $\bb_d$ which has entries $\bb_{d,j} = \int_\Omega \delta (\br-\br_d) u_j(\br) d\br$.

 The measurements that are collected at the detector are the values of the photon fluence produced as a result of different sources excited at various wavelengths. Let us denote by $\by$ the vector of measurements obtained as  
\begin{equation} \by = \sum_{k=1}^{\nsp} c_i\bE_i\bH\bmu + \bEta \qquad \bEta \sim \mathcal{N} (\bzero,\bW^{-2}) 
\label{eqn:measurement}
\end{equation}
where the rows of $\bH$ are constructed by the discretized representation of the integral  equation~\eqref{eqn:lippschwin} and $\bmu$ is a discrete representation of the shape of the  absorption perturbation $\chi(\br)$. The matrices $\bE_i = \bI_{N_s}\otimes \bI_{N_{ds}} \otimes \text{diag}(\varepsilon_k(\lambda))$ for $k=1,\dots,\nsp$ and $\varepsilon_k$ are the extinction coefficients for species $k$ that is a function of wavelength. Furthermore, $c_k$ are the concentration of the $k$-th chromophore species. The measurements are typically corrupted by noise which we model as a Gaussian distribution $\mathcal{N}(\bzero,\bW^{-2}$). The noise covariance $\bW$ is modeled as a diagonal matrix with diagonal entries $1/\sigma_m$. The relationship between the standard deviation $\sigma_m$ to the signal to noise ratio (SNR) is described in~\cite{larusson2011hyperspectral}.

\section{Fast solvers for the Born approximation}\label{sec:krylov}

In order to construct the measurement operator $\bH$ and thereby solving the inverse problem, we need the solution of the incident field $\bphii$ and $\bphid$ corresponding to multiple source and detector locations and multiple wavelengths. In this Section, we will discuss an efficient solver for the computing the sequence of parametrized systems using a novel recycling approach  based on Krylov subspaces. Krylov subspace methods are a popular class of algorithms for iterative solution to linear systems. Recycling Krylov subspaces in the context of parametric systems with smoothly varying coefficients has been previously considered in~\cite{chan1999galerkin,kilmer2006recycling,parks2006recycling}. %

\subsection{Recycling across wavelengths}
For simplicity of notation, we denote by \[ \sigma_j \define \frac{\nu\mu_a(\lambda_j)}{D(\lambda_j)} \qquad \text{and} \qquad \sigma_j' \define \frac{1}{2AD(\lambda_j)}\] the shifts corresponding to wavelengths $j=1,\dots,N_\lambda$,  and  by $\bx_j \define \bphii D(\lambda_j)$ so that equation~\eqref{eqn:discincident}
\begin{equation}\label{eqn:multipleshifted} 
(\bK + \sigma_j \bM + \sigma_j'\bR) \bx_j = \bb \qquad j = 1,\dots, N_\lambda \end{equation} 
where matrices $\bK,\bM$ and $\bR$ and vector $\bb$ are independent of the shift $\sigma$ and $\sigma'$. The solution to the adjoint field $\bphid$ can be written in a similar fashion with a different right hand side $\bb_d$. We first make the following change of variables $\bK \leftarrow \bM^{-1/2}\bK\bM^{-1/2}$,  $\bR \leftarrow \bM^{-1/2}\bR\bM^{-1/2}$ and transform the vectors as $\bx_j \leftarrow \bM^{-1/2}\bx$ and $\bb \leftarrow \bM^{-1/2}\bb$. This can be done efficiently by using a lumped mass matrix~\cite{hughes2012finite}. Denoted by $\krylov{\bK}{\bb}{n}$, the Krylov subspace of the matrix $\bK$ with starting vector $\bb$, is defined as 
\[ \krylov{\bK}{\bb}{n} \define \Span{\bb,\bK\bb,\dots,\bK^{n-1}\bb}\] 
Krylov subspaces enjoy an interesting property called shift invariance~\cite{simoncini2007recent}, i.e. \[ \krylov{\bK}{\bb}{n} = \krylov{\bK+ \sigma\bI}{\bb}{n}\]
 Several efficient methods exist for solving the system of equations~\eqref{eqn:multipleshifted} (if $\bR = 0$ corresponding to Neumann b.c.s), which solve for multiple shifts roughly at the cost of solving a single system. This is accomplished by generating a subspace that is independent of the shift and use the  shift-invariant property of Krylov subspaces (for a detailed review, see~\cite[Section 14.1]{simoncini2007recent} and references therein). However, the presence of a third matrix $\bR$ destroys the shift-invariant property of the Krylov subspace methods, unless $\krylov{\bK}{\bb}{n}$ is an invariant subspace of $\bR$. Even though the shift-invariant property does not hold, \textit{we can utilize information from the solution of the shifted system of equations $(\bK+\sigma_j\bI)\bx_j = \bb$ if the perturbation $\bR$ is not too large in magnitude}. This is the main idea behind the recycling approach to Krylov subspaces that we are proposing.

We construct a shift-invariant basis  $\krylov{\bK}{\bb}{n}$ by running $n$ steps of the Arnoldi algorithm and we extract the $k$ smallest Harmonic Ritz eigenvalues and eigenvectors to construct $\bU$ and $\bC$ both in $\mathbb{R}^{N\times k}$ such that 
\begin{equation}\label{eqn:gcr} \bK\bU = \bC \qquad \bC^T\bC = \bI
\end{equation}
By using the shift-invariant property of Krylov subspaces, we know that $\bU$ is also an approximately invariant subspace of $\bK+\sigma_j\bI$.  We now consider the extension of the solution to the parametric system $\bA_j \define \bK + \sigma_j \bI + \sigma_j'\bR$. When $\normtwo{\bR}$ is small relative to $\normtwo{\bK}$, this can be considered a perturbation to the system $\bK + \sigma_j\bI$ for which we have already generated an approximately invariant subspace. To compute a relation of the kind in Equation~\eqref{eqn:gcr} for the matrix $\bA_j$, we proceed as 
\[\bA_j \bU = \bC + \sigma_j \bU + \sigma_j'\bR \bU  \define \bC_j'\]
In the above equation, the most expensive step is computing $\bR_\bU \define \bR \bU$. The matrix $\bR$ only has non-zero entries corresponding to boundaries at which there is refraction index mismatch, in our application it is limited to the top and the bottom boundaries. As a consequence, $\bR$ is even more sparse compared to $\bK$. Moreover, the matrix $\bR\bU$ can be precomputed since it will be used across each shift. Compute the thin QR decomposition  which using MATLAB notation we represent as $[\bQ_j,\bY_j] = \text{qr}(\bA_j \bU, 0)$. The updated updated matrices can now be computed as  $\bU_j \define \bU\bY^{-1}_j$ and $\bC_j = \bQ_j$ and satisfy the relation 
\[\bA_j \bU_j = \bC_j \qquad \bC_j^T\bC_j = \bI. \]
Here and henceforth, $\bU_j$ is not computed explicitly, rather a solve using the upper triangular matrix $\bY_j$ is performed when it is necessary to form products with $\bU_j$. 

\begin{algorithm}[!ht]
\begin{algorithmic}[1]
\STATE\label{shiftinvariant} Generate a basis for Krylov subspace $\krylov{\bK}{\bb}{n}$ 
\[ \bK \bV_n = \bV_{n+1}\bar{\bT}_n \qquad \by_j \define \argmin \left\Vert \bV_{n+1}^T\bb - \left(\bar{\bT}_{n+1} + \sigma_j \begin{bmatrix}\bI \\ \bzero \end{bmatrix}\right) \by \right\Vert\]
\STATE Solve the eigenvalue problem $\bar{\bT}_n^T\bar{\bT}_n \bz = \theta \bT_n^T\bz $ and retain $k$ smallest eigenvalues in magnitude. Collect the eigenvectors into a matrix $\bZ_k$ and the corresponding eigenvalues as $\boldsymbol\Theta_k$.
\STATE Compute $\bC  =   \bV_{n+1}\bar{\bT}_n\bZ_k\boldsymbol\Theta_k$ and $\bU_k = \bV_{n+1}\bar{\bT}_n\bZ_k$
\STATE  Compute the thin QR $[\bC,\bY] = \text{qr}(\bC,0)$ and set $\bU = \bU\bY^{-1}$
 \COMMENT {Shift-invariant deflation subspace}
\STATE Compute the initial solution $\tilde{\bx}_{0,j} = \bV_n\by_j$ for $j=1,\dots, N_\lambda$ 
\FOR {$j = 1,\dots,N_\lambda$}
\STATE  $\bC_j' = \bC + \sigma_j \bU + \sigma_j'\bR\bU $  \COMMENT{Compute new deflation subspace}
\STATE $[\bC_j, \bY_j] = \text{qr}(\bC_j',0)$ and set $\bU_j = \bU\bY_j^{-1}$.
\STATE Compute $\br_j = \bb - \bA_j\bx_{0,j}$ and $\bx_{-1,j} \define \tilde{\bx}_{0,j} + \bU_j\bC_j^T\br_j$ \COMMENT {Initial deflation}
\STATE Generate $\bV_{m-k+1}^{(j)}$ and $\bar{\bT}_{m-k}^{(j)}$ by applying $m-k$ steps of the Arnoldi algorithm using matrix $(\bI - \bC_j\bC_j^T)\bA_j$ applied to the initial deflated vector $(\bI - \bC_j\bC_j^T)\br_j$

\STATE Solve the least-squares system of equations~\eqref{eqn:minres} with $\beta_j = \normtwo{(\bI - \bC_j\bC_j^T)\br_j}$
\[ [\bQ,\bZ] = \text{qr}(\bar{\bT}_{m-k}^{(j)},0) \quad \by_{2,j} = \bZ^{-1}(\bQ^T\beta_j \be_1) \quad \by_{1,j} = - \bF_k^{(j)}\by_{2,j}  \] 
to generate approximate solution $\bx_{m,j} = \bx_{-1,j}  + \bU_j\by_{1,j} + \bV_{m-k}^{(j)}\by_{2,j}$
\ENDFOR
\end{algorithmic}
\caption{Augmented GMRES for parametric systems}
\label{alg:auggmres}
\end{algorithm}

We adopt the same recycling strategy as~\cite{parks2006recycling}. However, where we differ from this approach is the way we construct and update the approximate invariant subspace, as we now explain. The strategy in~\cite{parks2006recycling} is to first generate an approximate a recycling basis for the matrix $\bA_1$ and then update the recycling subspace of perturbed matrices $\bA_j = \bA_{j-1} + \Delta \bA_j$ after the convergence of the augmented Krylov solver for each of the system.   In the context of our problem, note that 
\[\Delta \bA_j \define (\sigma_j-\sigma_{j-1})\bI + (\sigma_j' -\sigma_{j-1}')\bR .\]
 Since the matrices are known a priori and only the shifts are varying, we are able to do something different: once the basis $\bU$ and $\bC$ are known, our approach to update $\bU_j$ and $\bC_j$ can be easily performed independently and can be parallelized in a straightforward manner.   We will discuss the steps and cost for obtaining
 $\bC_j$ at the end of this section.  First, let us assume that $\bU_j, \bC_j$ are available, and discuss our solution technique.

To begin, we assume we have computed  
$\tilde{\bx}_{0,j}$ iteratively and simultaneously for all $j$ using the shift-invariant property as the approximate
solution to 
$(\bK + \sigma_j \bI) \bx_j= \bb$.  We use $\tilde{\bx}_{0,j}$ as a first estimate of the solution to  $(\bK + \sigma_j \bI + \sigma_j' \bR) \bx = \bb$ obtained by exploiting the shift-invariance of the Krylov subspaces (see Step~\ref{shiftinvariant} in Algorithm~\ref{alg:auggmres}).   The corresponding initial residual is denoted $\br_{0,j} := \bb - \bA_j \tilde{\bx}_{0,j}$. Because we expect this solution to need augmentation, we next search for a 
better estimate of the form $\bx_{-1,j} = \tilde{\bx}_{0,j} + \bU_j \bz$.  We choose $\bz$ such that
\[ \bz = \argmin \| \bb - \bA_j(\tilde{\bx}_{0,j} + \bU_j \bz ) \| = \argmin \| \br_{0,j} - \bC_j \bz \|    \]
 For this choice, we get a new solution estimate $\bx_{-1,j} = \tilde{\bx}_{0,j} + \bU_j \bC_j^T (\bb - \bA_j\tilde{\bx}_{0,j} )$ that gives the residual $\br_{-1,j} = (\bI - \bC_j \bC_j^T) \br_{0,j}$.


The approximate solution $\bx_{m,j}$ for the parametric system $j$ is obtained by searching in the augmented affine subspace 
\begin{equation}
\bx_{m,j} \in \bx_{-1,j} + \Span{\bU_j} \oplus \krylov{(\bI-\bC_j\bC^T_j)\bA_j}{(\bI-\bC_j\bC^T_j)\br_{0,j}}{m-k}
\end{equation}
that is, by searching for solutions in the deflated subspace $\krylov{(\bI-\bC_j\bC^T_j)\bA_j}{(\bI-\bC_j\bC^T_j)\br_{0,j}}{m-k}$ obtained by applying $m-k$ steps of the  Arnoldi algorithm with the deflated matrix $(\bI-\bC_j\bC^T_j)\bA_j$ to the starting vector $(\bI-\bC_j\bC^T_j)\bb$ and augmented with the subspace $\bU_j$. The initial guess $\bx_{-1,j}$ is obtained as $ \bU_j\bC^T_j(\bb - \bA\tilde{\bx}_{0,j})$ and the initial residual is $\br_{0,j} = \bC_j\bC^T_j\bb$, (where $\tilde{\bx}_{0,j}$ is the approximation solution obtained by using the shift-invariant property to solve $(\bK + \sigma_j\bI)\bx_j = \bb$).

The following discussion closely mirrors the presentation of recycling in \cite{parks2006recycling,wang2007large,kilmer2006recycling,mello2010recycling}.   
The Arnoldi algorithm on the deflated problem yields the matrix relationship

\[ (\bI - \bC_j\bC^T_j)\bA_j \bV_{m-k}^{(j)} = \bV_{m-k+1}^{(j)}\bar{\bT}_{m-k}^{(j)}\]
where the superscripts indicate the system index $j$. The above equation can be rewritten as a  modified Arnoldi relationship by defining  $\bF_k^{(j)} \define \bC^T_j\bA_j\bV_{m-k}^{(j)}$ and reorganizing as
\begin{equation}
\bA [\bU_j, \bV_{m-k}^{(j)}] = [\bC_j, \bV_{m-k+1}^{(j)}] \begin{bmatrix}\bI_k & \bF_k^{(j)} \\ \bzero & \bar{\bT}_{m-k}^{(j)} \end{bmatrix}
\end{equation}
Now, $\bV_{m-k}^{(j)}$ forms a basis for the subspace $\krylov{(\bI-\bC_j\bC^T_j)\bA_j}{(\bI-\bC_j\bC^T_j)\br_{0,j}}{m-k}$ and now we search for solutions of the form $\bx_{m,j} \in \bx_{-1,j} + \Span{\bU_j} \oplus \Span{\bV_{m-k}^{(j)}}$ and can be written as 
\[ \bx_{m,j} = \bx_{-1,j} + \bU_j\by_{1,j} + \bV_{m-k}^{(j)} \by_{2,j}\]
The solution to the coefficients $\by_{1,j}$ and $\by_{2,j}$ are obtained by minimizing the residual which results in the following least squares problem
\begin{equation}\label{eqn:minres}  \min_{\by_{1,j},\by_{2,j}} \left\lVert \begin{bmatrix} \bzero \\ \beta_j \be_1 \end{bmatrix} - \begin{bmatrix}\bI_k & \bF_k^{(j)} \\ \bzero & \bar{\bT}_{m-k}^{(j)} \end{bmatrix} \begin{bmatrix} \by_{1,j} \\ \by_{2,j} \end{bmatrix}\right\rVert  \end{equation}
where $\beta_j = \normtwo{(\bI-\bC_j\bC^T_j)\br_{0,j}}$. Since the number of iterations are expected to be small, we store the vectors $\bV_{m-k}^{(j)}$ and solve the least squares problem in equation~\eqref{eqn:minres} directly. If the number of iterations are expected to be large, we do not need the solution to $\by_{1,j}$ and $\by_{2,j}$ explicitly, only the products $\bV_{m-k}^{(j)}\by_{2,j}$ and $\bF_k^{(j)} \by_{2,j}$ which can be obtained using short-term recurrence relation similar to recycled MINRES~\cite{wang2007large} (assuming the matrices $\bK$, $\bM$ and $\bR$ are symmetric).

\subsection{Computational and storage costs}
We now discuss the computational and storage costs involving Algorithm~\ref{alg:auggmres} and the overhead induced by the augmented approach for constructing the augmented basis $\bU$. We note that the loop can be executed in parallel because no information is shared across the solves except for the initial choice of $\bU$. We focus only on costs that are linear in the size of the matrix, i.e. costs of the form $\bigO(N)$ since the costs involving smaller matrices are negligible. The augmented approach in Algorithm~\ref{alg:auggmres} requires additional storage of $2Nk$ for the matrices $\bU$ and $\bC$.  For the pre-computation of the new bases $\bU_j$ and $\bC_j$ the major cost is the QR factorization which is $\bigO(Nk^2)$ for each wavelength. To accelerate this computation, we use a more efficient approach at the possible expense of some accuracy.  We first compute
 \begin{align*}(\bC_j')^T\bC_j' = & \quad \bI + \sigma_j ( \bC^T\bU + \bU^T\bC)+\sigma_j'( \bC^T\bR_\bU+ \bR_\bU^T\bC)  \\
& \quad  + \sigma_j\sigma_j' (\bR_\bU^T\bU + \bU^T\bR_\bU) + \sigma_j^2 \bU^T\bU+  (\sigma_j')^2 \bR_\bU^T\bR_\bU .
\end{align*}
Next, the small $k\times k$ matrices such as $\bU^T\bC$ are precomputed and stored since they are independent of the shifts $\sigma_j$ and $\sigma_j'$. Then the Cholesky factorization of $(\bC_j')^T\bC_j'$ is computed and  $\bC_j$ is obtained as $\bC_j = \bC_j'\bF^{-1}$. Note as before that the inverse $\bF^{-1}$ is not computed explicitly. Because of this pre-computation, the additional cost per wavelength is now only $\bigO(k^3)$. Suppose the algorithm converges in $m$ iterations, then the algorithm requires $m-k$ additional matrix-vector products and $\bigO(Nm^2 + mkN )$ other floating point operations. For an efficient algorithm, the overhead costs must be offset by the gains obtained by decreasing the number of iterations. Numerical experiments performed in Section~\ref{sec:results} demonstrate that the reduction in the number of matrix-vector products due to deflation offsets the additional computational cost due to pre-computation and re-orthogonalization.

\section{Fast compression of $\bH$}\label{sec:compress}
After computing the incident field $\bphii$ and the adjoint field $\bphid$ the next step is to compute the measurement operator $\bH$. The construction of $\bH$ is described in Equation~\eqref{eqn:measurement}. However, recall that the operator $\bH$ is expensive to store and compute, and therefore our goal is to produce an approximate factorization of $\bH \approx \bU\bV^T$. Furthermore, an optimal low-rank compression using SVD scales as $\bigO \left( \min \{M,N\}^2, \max\{M,N\} \right)$ which is prohibitively expensive. Here we present an algorithm that avoids computing $\bH$ in its entirety, but computes and compresses sub-blocks and then combines the compressed sub-blocks in a recursive manner.

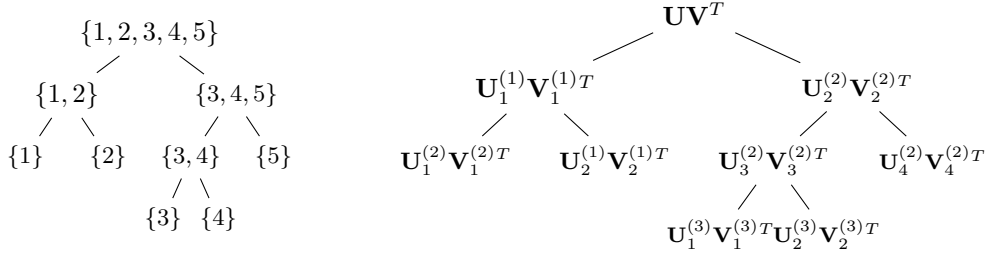
\begin{figure}
\centering
\begin{subfigure}{0.4\textwidth}\centering
\begin{tikzpicture}[level/.style={sibling distance=40mm/#1},scale=0.55]
    \tikzstyle{ann} = [draw=none,fill=none]
\node [ann] (z){$\{1,2,3,4,5 \}$}
  child {node [ann] (a) {$\{1,2\}$}
    child {node [ann]  {\small{$\{ 1\}$}}
    }
    child {node [ann]  {\small{$ \{ 2\} $}}
    }
  }
child {node [ann] (j) { \small{ $\{ 3,4,5\}$}}
    child {
    	node [ann]  {\small{$\{3,4\}$}}
    	child { 
    		node [ann] {\small{$\{3\}$} }
    	}
    	child{
    		node [ann] {\small{$\{4\}$} }
    	}
    }
  child {
  	node [ann] (l) {\small{$\{5\}$}}
  }
};
\end{tikzpicture}
\end{subfigure}
\begin{subfigure}{0.4\textwidth}
\begin{tikzpicture}[level/.style={sibling distance=65mm/#1},scale=0.65]
    \tikzstyle{ann} = [draw=none,fill=none]
\node [ann] (z){$ \bU\bV^T$}
  child {node [ann] (a) {$\bU_1^{(1)}\bV_1^{(1)}{}^T$}
    child {node [ann]  {\small{$\bU_1^{(2)}\bV_1^{(2)}{}^T$}}
    }
    child {node [ann]  {\small{$ \bU_2^{(1)}\bV_2^{(1)}{}^T $}}
    }
  }
child {node [ann] (j) { \small{ $\bU_2^{(2)}\bV_2^{(2)}{}^T$}}
    child {
    	node [ann]  {\small{$\bU_3^{(2)}\bV_3^{(2)}{}^T$}}
    	child { 
    		node [ann] {\footnotesize{$\bU_1^{(3)}\bV_1^{(3)}{}^T$} }
    	}
    	child{
    		node [ann] {\footnotesize{$\bU_2^{(3)}\bV_2^{(3)}{}^T$} }
    	}
    }
  child {
  	node [ann] (l) {\footnotesize{$\bU_4^{(2)}\bV_4^{(2)}{}^T$}}
  }
};
\end{tikzpicture}
\end{subfigure}

\caption{The size of the matrices at the leaf level are $ N_{ds}N_\lambda \times N$ and $N_s = 5$. The superscripts denote the level of the tree, where as the subscripts denote the order at a given level (which can be different than the global ordering of the sources). At the leaves the matrices are approximated by low-rank factors and then agglomerated recursively based on the tree. Here we consider a tree corresponding to $5$ sources. A possible source configuration that produces the tree can be obtained in Figure~\ref{fig:domain}. }
\label{fig:part}

\end{figure}

\subsection{Outline} The rows of $\bH$ can be partitioned as 

\[ \bH^T = [\bH_1^T,\dots,\bH_{N_s}^T] \]
where each block $\bH_i$ is of size $N_{ds}N_\lambda \times N$. Each row of $\bH$ represents a discretized version of the integral described in Equation~\eqref{eqn:lippschwin} that combines the incident field and the adjoint field. Instead of compressing $\bH$ by using low-rank factorization techniques such as truncated SVD, the idea is to compress each block $\bH_i \approx \bU_i\bV_i^T$ locally and then combine the factorizations in a recursive fashion. The full algorithm is provided in Algorithm~\ref{alg:recursive} and is illustrated for $N_s = 4$ in Figure~\ref{fig:part}. The low-rank approximation of the blocks $\bH_i$ can be accomplished either using Randomized SVD or partially pivoted ACA algorithms described in Subsection~\ref{sec:compression}. Then a scheme for ordering the blocks $\bH_i$ is presented in Subsection~\ref{sec:part} that uses a spatial bisection tree to order the source locations by their spatial proximity. The agglomeration of the low-rank factors is accomplished by recursion using this tree structure. In Subsection~\ref{sec:costs} we analyze the computational costs of this recursive compression scheme and conclude with an error analysis in Subsection~\ref{sec:acc}. 

\subsection{Low Rank representation}\label{sec:compression}

As was mentioned earlier, we require a strategy to compute a low-rank factorization of the sub-blocks $\bH_j^T$ for $j=1,\dots,N_s$. Consider the blocks  $ \bH_j $ of size $m \times n$ where $m = N_{ds}N_\lambda$ and $n = N$. Considering only one source, we seek a low-rank approximation of the form $\bH \approx \bU \bV^T$ and the number of columns of $\bU$ and $\bV$ is denoted by $r$. It is well known that the best low-rank approximation to rank $r$ is obtained by truncating the SVD to rank $r$. In this case, we have $\| \bH - \bU_r\Sigma_r\bV_r^T \|_2  = \sigma_{r+1} $.  However, computing the SVD is expensive since it requires $\bigO(nm^2)$ operations assuming $n \geq m$. The advantage of the low-rank representation is that the cost of storing the decomposition and computing matrix-vector products are both given by $\bigO(r(m+n))$ instead of $\bigO(mn)$. When $r \ll \min\{m,n\}$, this reduction can represent significant savings.

In the Appendix~\ref{sec:lowrank} we describe two approximate methods that compute a low-rank representation but have a lower computational cost asymptotically than the SVD, namely randomized SVD (RandSVD) and partially pivoted Adaptive Cross Approximation (ppACA). Here, we only summarize the resulting computational costs. For a matrix of dimensions $m\times n$ that has a rank $r$ the cost can be expressed as 
\begin{equation}
\mu_\text{Comp}(m,n;r) \quad = \quad \left\{ \begin{array}{cc} C_1rmn + C_2 r^2(m+n) & \text{RandSVD}\\ C_3r^2(m+n)&  \text{ppACA}\end{array}\right.
\label{eqn:svdcosts}
\end{equation}
where constants $C_1$, $C_2$ and $C_3$ are assumed to be known and provided in the literature (see Appendix~\ref{sec:lowrank}). Note that the costs of the low-rank approximation both using RandSVD and ppACA are asymptotically smaller than the cost of an SVD which scales as $\bigO(nm^2)$ assuming $n > m$.

\subsection{Agglomerating low ranks}\label{sec:agglomeration}

Having produced low-rank approximations to the sub-blocks $\bH_j$ for $j = 1,\dots,N_s$ we now consider the problem of agglomerating low-rank factors to produce a global low-rank factorization. Here we consider only two sources, i.e. $N_s = 2$ and as before, $\bH^T = [\bH_1^T, \bH_2^T]$ of size $m \times n$ where $m =  N_{ds}N_\lambda$ and $n = N$.  Suppose we have the low-rank factors $\bH_1 \approx \bU_1 \bV_1^T$ and $\bH_2 \approx \bU_2 \bV_2^T$ each with rank $r$ which have been compressed according to some predetermined tolerance $\varepsilon$. The low-rank factors can be then combined as 
\[\bH = \begin{bmatrix} - \bH_1 -  \\ - \bH_2 - \end{bmatrix} \approx \begin{bmatrix}  \bU_1\bV_1^T \\  \bU_2\bV_2^T \end{bmatrix} = \begin{bmatrix}\bU_1 & \\ & \bU_2 \end{bmatrix} \begin{bmatrix} - \bV_1^T - \\ - \bV_2^T -  \end{bmatrix}  .\]
Now the leftmost matrix is $m \times r$ and has independent columns by construction.  
The rightmost $r \times n$ matrix, however, may have a rank smaller yet than $\min{r,n}$.  So we compute 
$\begin{bmatrix} - \bV_1^T-  \\ -  \bV_2^T -  \end{bmatrix} = \bU_\bV \bV^T$ where $\bU_\bV, \bV$ have
$r' < r$ columns,  and we set
\[  \begin{bmatrix}\bU_1 & \\ & \bU_2 \end{bmatrix} \bU_\bV \bV^T = \bU\bV^T, \qquad \bU \in \mathbb{R}^{2N_{ds}N_{\lambda} \times r'} , \bV \in \mathbb{R}^{N \times r'} .\]

We will need to form $\bU$ explicitly, and this requires computing a rank-r' approximation to the stacked $\bV_i$ matrix.  The cost must therefore include the matrix-matrix product (using the block diagonal structure of the leftmost matrix in the product) and the cost of this additional factorization, using the same algorithm as was used
to approximately factor the $\bH_i$.  
The estimate of the cost of agglomeration of sub-blocks is therefore,
\begin{equation}
\mu_\text{Agg}(m,n;r',r) \quad = \quad  
\left\{ \begin{array}{cc} C'rr'(m+n) & \text{RandSVD}\\ C''(r')^2(m+n)&  \text{ppACA}\end{array}\right.
\end{equation}
where the constant $C'$ and $C''$ are assumed to be known (see Appendix~\ref{sec:lowrank}). For simplicity of analysis we will only use RandSVD for agglomerating low-rank representations.

\subsection{Recursive SVD}\label{sec:part}
We have, thus far, described a procedure to compute low-rank factors $\bH_j \approx \bU_j\bV_j^T$ corresponding to source indexes $j=1,\dots,N_s$ and shown how to agglomerate low-rank factors when the number of sources are $2$. However, as we shall show, this procedure can be implemented recursively. The complete algorithm for computing and compressing $\bH$ is summarized in Algorithm~\ref{alg:recursive}.

Let $I  = \{1,\dots,N_s\}$  denote an index set and let $|I|$ denote the cardinality of the set $I$.  Let the index sets $I$ be partitioned into binary trees denoted by $T_I$  respectively. For each $\tau \in T_I$, we denote the set of its sons by $S(\tau) \subset T_I$ and the leaves of the tree are denoted by $\mathcal{L}(T_I)$. The tree $T_I$ can be used to recursively spatially partition the domain so that the sources are ordered such that recursively combining low-rank factors from sub-blocks corresponds to combining nearby sources since the corresponding Green's functions are known to be highly compressible~\cite{chaillat2012faims}. As a result, by this ordering, we hope to gain a larger compression at each level in the tree. Given the tree $T_I$ we can recursively compress the sub-blocks using Algorithm~\ref{alg:recursive} which is initialized using the root of the tree, which corresponds to the index set $I$. The partitioning of the matrix $\bH$ into sub-blocks and their recursive compression for $N_s = 4$ using the Figure~\ref{fig:part}. The construction of the tree is described in Algorithm~\ref{alg:recursive}. 

\begin{algorithm}[!ht] \label{alg:rr}
\begin{algorithmic}[1] 
\IF {$S(\tau) \neq \emptyset$}
\STATE $\mathcal{U} := \emptyset$, $\mathcal{V} := \emptyset$
\FORALL {$\tau' \in S(\tau ) $}
\STATE  $\bU_{\tau'},\bV_{\tau'}  = \text{RecursiveLowRank}(\tau', \varepsilon)$
\STATE $\mathcal{U}.\text{append} ( \bU)$, $\mathcal{V}.\text{append} ( \bV)$
\ENDFOR
\ELSE 
\FORALL {$\tau' \in \mathcal{L}(T_I ) $}
\STATE Compute the sub-block $\bH_\tau'$ using Equation~\eqref{eqn:lippschwin} \\
\COMMENT {// Use fast solvers developed in Section~\ref{sec:krylov} for incident  $\bphii$ and adjoint field $\bphid$. }
\STATE $[\bU,\bV] = \text{LowRank}(\bH|_{\tau'}, \varepsilon)$  \\
\COMMENT {//Compute low-rank at leaf level such that $\norm{\bH_{\tau'}-\bU\bV^T} \leq 
\varepsilon \norm{\bH_{\tau'}}$. see Section~\ref{sec:compression}}
\STATE $\mathcal{U}.\text{append} ( \bU)$, $\mathcal{V}.\text{append} ( \bV)$
\ENDFOR
\ENDIF
\STATE $[ \bU,\bV]=$ Agglomerate($\mathcal{U},\mathcal{V}$, $\varepsilon$)
\COMMENT {// Agglomerate low-rank factors, see Section~\ref{sec:agglomeration}}

\RETURN $\bU, \bV$ such that $\bH \approx\bU\bV^T$
\end{algorithmic}
\caption{RecursiveLowRank$(\tau, \varepsilon)$}
\label{alg:recursive}
\end{algorithm}

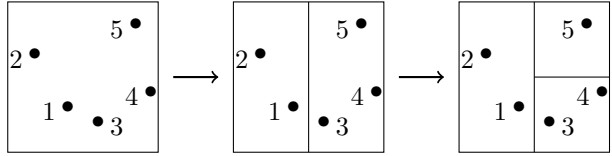
\begin{figure}
 \centering
 \begin{subfigure}{0.75\textwidth}\centering
 \begin{tikzpicture}[scale = 0.2]
 \def\l{10}
 \def\s{ 5}
\draw (0,0) rectangle (\l,\l);
\foreach \Point [count = \xi] in {(0.4*\l,0.3*\l), (0.18*\l,0.65*\l), (0.6*\l,0.2*\l), (0.95*\l,0.4*\l), (0.85*\l, 0.85*\l)}{
 \ifnum\xi=3  
 	\node [label = {[xshift=0.025*\l cm,yshift=-0.05*\l cm]$\xi$}] at \Point {\textbullet};
 \else  	
    	\node [label = {[xshift=-0.025*\l cm,yshift=-0.05*\l cm]$\xi$}] at \Point {\textbullet};    
    \fi
}
\draw [->, thick] (\l +0.2*\s, .5*\l) -- (\l +0.8*\s,.5*\l); 
\draw (\l+\s,0) rectangle (2*\l +\s , \l);
\draw (1.5*\l + \s,0) -- (1.5*\l + \s,\l);
\foreach \Point [count = \xi] in {(\l+\s + 0.4*\l,0.3*\l), (\l+\s + 0.18*\l,0.65*\l), (\l+\s + 0.6*\l,0.2*\l), (\l+\s + 0.95*\l,0.4*\l),  (\l+\s + .85*\l, 0.85*\l)}     {
    \ifnum\xi=3  
 	\node [label = {[xshift=0.025*\l cm,yshift=-0.05*\l cm]$\xi$}] at \Point {\textbullet};
 \else  	
    	\node [label = {[xshift=-0.025*\l cm,yshift=-0.05*\l cm]$\xi$}] at \Point {\textbullet};    
    \fi
}
\draw [->, thick] (2*\l + \s +0.2*\s, .5*\l) -- (2*\l + \s +0.8*\s,.5*\l); 
\draw (2*\l + 2*\s,0) rectangle (3*\l + 2*\s ,\l);
\draw (2.5*\l + 2*\s,0) -- (2.5*\l + 2*\s,\l);
\draw (2.5*\l + 2*\s, .5*\l) -- (3*\l + 2*\s,0.5*\l);

\foreach \Point  [count = \xi] in {(2*\l+2*\s + 0.4*\l,0.3*\l), (2*\l+2*\s + 0.18*\l,0.65*\l), (2*\l+2*\s + 0.6*\l,0.2*\l), (2*\l+2*\s + 0.95*\l,0.4*\l),  (2*\l+2*\s + .85*\l, 0.85*\l)}{
   \ifnum\xi=3  
 	\node [label = {[xshift=0.025*\l cm,yshift=-0.05*\l cm]$\xi$}] at \Point {\textbullet};
 \else  	
    	\node [label = {[xshift=-0.025*\l cm,yshift=-0.05*\l cm]$\xi$}] at \Point {\textbullet};    
    \fi	
}

\end{tikzpicture}
\end{subfigure}

\caption{ Top view of the bounding box for a possible source configuration. Locations of $5$ different sources that can be subdivided recursively using Algorithm~\ref{alg:geom}. The  resulting tree $T_J$ that is returned by the algorithm, when initialized by the index set $J = \{1,\dots,5\}$ is provided in Figure~\ref{fig:part}.  }
\label{fig:domain}
\end{figure}

We make the assumption that the locations of the source, detectors and the support of the perturbation are non-overlapping and well separated.  For concreteness, let us assume that the breast is placed between parallel plates and is enclosed by a cuboid of dimensions $[-L_x,L_x] \times [-L_y,L_y] \times [0,L_z]$. Furthermore, let us assume that the sources are located on the top plate $z=L_z$ and the detectors are located at $z=0$ and the detectors for any given source are roughly co-axial. The construction of the tree is performed as follows, and a simple example illustrating this construction is provided in Figure~\ref{fig:domain}.  Let $J$ be the index set corresponding to source locations. The tree $T_J$ is constructed recursively by geometric bisection applied on a 2D domain.  Given the initial bounding box containing all the source (in this case, a box of size $[-L_x,L_x] \times [L_y, L_y]$), the domain is split in a direction perpendicular to the direction of maximum expanse  and the sources are split between the newly created split domains. If the number of points in each domain are greater than $2$, the  procedure is computed recursively until the number of points in each leaf is no greater than $2$. If the points are uniformly distributed, the number of levels $L = \log_2(N_s)$, where $N_s$ is the number of sources. The algorithm is summarized in Algorithm~\ref{alg:geom}. The algorithm is initialized using $\tau = J$, $\alpha = [-L_x,L_x]$ and $[-L_y,L_y]$.

\begin{algorithm}[!h]
\begin{algorithmic}[1]
\STATE $j_\text{max} :=\arg\max\{ \beta_j- \alpha_j : j \in \{ 1,\dots,d\} \}$
\STATE  $\gamma := (\alpha_{j_\text{max}} + \beta_{j_\text{max}})/2$ \COMMENT{// Split cluster in chosen direction }
\STATE $\tau_1 := \emptyset$, $\tau_2 := \emptyset$ 
\FOR {$i \in \tau$}
\IF {$e_{j_\text{max}}^T\bx_i  \leq \gamma $ }
\STATE  $\tau_1 = \tau_1 \bigcup \{i\} $
\ELSE 
\STATE  $\tau_2 = \tau_2 \bigcup \{i\} $
\ENDIF
\ENDFOR
\STATE Define $\beta' := \beta$, $\beta_{j_\text{max}}' := \gamma$,   and $\alpha' := \alpha$, $\alpha_{j_\text{max}}' := \gamma$
\STATE Split$(\tau_1,\alpha,\beta') $ and Split$(\tau_2,\alpha',\beta)$ \COMMENT {// Split only if $|\tau_i| > 2$ for $i=1,2$.}
\end{algorithmic}
\caption{Split$(\tau, \alpha, \beta)$ // Geometric bisection to construct tree}
\label{alg:geom}
\end{algorithm}

\subsection{Computational cost}\label{sec:costs}
In this section we present a detailed analysis of the computational cost of the recursive SVD described in Algorithm~\ref{alg:recursive}. Our analysis is in the same spirit as the one provided in~\cite{hackbusch2014new}. However, while they used the full SVD algorithm to produce the optimal truncation at each level, we will consider the use of approximate low-rank factorizations that have better computational costs. As we shall see, this analysis and results are more sharp than~\cite{hackbusch2014new}  
since we account for the cost depending on rank of the sub-blocks at each level of the tree.

The starting point of our analysis  is splitting the costs into a contribution that comes from the preparation at the leaf level that is computed using the low-rank factorization techniques listed in Section~\ref{sec:compression} and the contribution that comes from agglomerating low-rank sub-blocks.  The number of sub-blocks to agglomerate become fewer the higher up we are in the tree, but the size of the matrices to be agglomerated increases. The total cost is therefore equal to 

\begin{equation}\label{eqn:genericcosts}
\text{Cost } = \sum_{b \in \mathcal{L}(T_{I})} \mu_{Comp}(m^{(L)},n^{(L)};r_L) + \sum_{\ell=0}^{L-1}\sum_{b \in T^{(l)}_{I}  } \mu_\text{Agg}\left(m^{(\ell)},n^{(\ell)};r_{L-\ell-1},r_{L-\ell}\right) \text{ flops}.
\end{equation}

We assume that $N_s$ is a power of $2$ and each partition has the same size, so that the size of the partition at the leaf level is $N_\lambda N_{ds} \times N$, and $N_s = 2^L$ where $L$ is the depth of the tree $T_{I}$. Furthermore, we denote by $r_{\ell}$, the maximum block rank of any partition at level $\ell$.  The cost of recursive SVD is then 
\begin{align*}
\text{Cost } = & \quad 2^L\mu_\text{Comp}\left(\frac{M}{2^L},N; r_L\right) & +& \quad  \sum_{\ell=0}^{L-1}2^{L-\ell-1}\mu_\text{Agg}\left(2^{\ell -L}M,N;r_{L-\ell - 1},r_{L-\ell}\right) \\ 
= & \quad C_12^Lr_L\frac{M}{2^L}N +  2^LC_2r_L^2\left(\frac{M}{2^L} + N\right) & +&  \quad \sum_{\ell=0}^{L-1}C'2^{L-\ell-1}r_{L-\ell-1}r_{L-\ell}(2^{\ell-L} M + N) \\ 
= & \quad  C_1r_LM{N} +  C_2r_L^2\left(M+2^LN\right) & + &\quad  \sum_{\ell=0}^{L-1}C'2^{L-\ell-1} 2\delta^{L-\ell-1}_{L-\ell}r_{L-\ell}^2(2^{\ell-L} M + N) \\
= & \quad C_1r_LM{N} +  C_2r_L^2\left(M+2^LN\right)  &+ & \quad 2^{L-1}C'r_{L}^2\sum_{\ell=0}^{L-1}(M/2^L + 2^{-\ell}N)2\delta^{L-\ell-1}_{L-\ell} \prod_{k=L-\ell}^{L-1}\left(2\delta^{k}_{k+1}\right)^2
\end{align*}
Here, we assume that $r_{\ell} = 2\delta^{\ell}_{\ell+1}r_{\ell+1}$, that is, the ratio of the ranks at level $\ell$ that is obtained by agglomeration of two sub-blocks each of rank at most $r_{\ell+1}$ is a factor of $\delta^{\ell}_{\ell+1}$, which is dependent on the particular level under consideration. It is easy to see that $0 < \delta^{\ell}_{\ell+1} \leq 1$ for all $\ell=0,\dots,L-1$. If we further make the assumption that $ \max_{0\leq \ell < L} \delta^{\ell}_{\ell+1} = \delta$ is independent of the level $\ell$, then we have following result that the total computational cost of recursive low-rank truncation can be further simplified to 
\begin{equation}
C_1r_LM{N} +  C_2r_L^2\left(M+2^LN\right) +   2^{L}C'r_{L}^2\delta\left(\frac{1-(\sqrt{2}\delta)^{2L-1}}{1-\sqrt{2}\delta}\frac{M}{2^L} +  \frac{1-(2\delta)^{2L-1}}{1-2\delta}{N}\right) .
\end{equation}
We define the quantities $f_1$ and $f_2$ as follows:
\begin{align*}
f_1 \define  \quad \sum_{\ell=0}^{L-1}2^{-\ell}\delta^{L-\ell-1}_{L-\ell}\prod_{k=L-\ell}^{L-1} \left(2\delta^{k}_{k+1}\right) \quad   &  \leq     \quad \delta \frac{1-(\sqrt{2}\delta)^{2L-1}}{1-\sqrt{2}\delta} \\
f_2 \define  \quad  \sum_{\ell=0}^{L-1}\delta^{L-\ell-1}_{L-\ell}\prod_{k=L-\ell}^{L-1} \left(2\delta^{k}_{k+1}\right) \quad  &  \leq    \quad \delta  \frac{1-(2\delta)^{2L-1}}{1-2\delta}.
\end{align*}
Here, we assume that $\delta \neq 1/\sqrt{2}, 1/2$. If this were the case, the appropriate sums would simplify to equal $L$. In the worst case, when there is no compression at higher levels $\delta = 1$ and then $f_1 \sim 2^L$ and $f_2 \sim 4^L$.

\begin{table}[!ht]\centering
\renewcommand{\arraystretch}{1.5}
\begin{tabular}{|c|c|c|c|c|} \hline
\multirow{2}{*} {Method} & \multicolumn{2}{|c|}{ Leaf computation} & \multicolumn{2}{|c|}{Tree computation} \\ \cline{2-5} 
           &RandSVD      & ppACA& Average & Worst \\ \hline
Direct & $C_1RMN + C_2 R^2 (M+N)$ & $C_3R^2(M+N)$ &  - &  - \\ \hline
Recursive &   $ C_1r_LMN + C_2r_L^2 (M+ 2^LN)$  & $C_3r_L^2(M+2^LN)$& $Cr_L^2(Mf_2 +2^LNf_1) $ & $r_L^2C4^L(M +N) $ \\ \hline
\end{tabular}
\caption{Summary of computational costs of the recursive SVD algorithm. Here $M=N_sN_{ds}N_\lambda$ is the number of measurements and $N$ is the grid size. The storage costs are $\bigO(R(M+N))$ where $R$ is the global rank of the low-rank factorization. } 
\label{tab:compcosts}
\end{table}

 We now compare the computational costs between different methods that we have outlined in this section. It can be readily seen that the cost of the recursive factorizations scale asymptotically better than the SVD which scales as $\bigO(NM^2)$ assuming $M \leq N$. We denote by `Direct', the low-rank algorithms described in Section~\ref{sec:compression}.  Considering only the costs that are of the order $\bigO(MN)$, it can be readily observed that simply  the recursive SVD methods have better scaling than the `Direct' method which scales as $\bigO(NM^2)$, if the rank at the leaves are smaller than the overall rank of the matrix $r_L < R$. By comparing the worst case costs for the recursive SVD methods, we can see that that the terms linear in $M$ and $N$ are comparable if $r_L \sim R/2^L$. If no compression is observed at any levels, including the leaves, then our algorithm performs poorly since we are needlessly computing a large number of ``low-rank'' factorizations at all the levels in the tree. However, there are still a couple of benefits of using the recursive SVD approach.  First,  for the range of parameters we are interested in exploring, storing the entire matrix $\bH$ could cost $\sim 200 $ GB which may be completely infeasible to store and later compress. By contrast the strategy in Algorithm~\ref{alg:recursive} does not require storage of $\bH$ in its entirety but computes and compresses sub-blocks of $\bH$ on-the-fly and therefore, has favorable storage costs. Second, our algorithm provides more locality in the calculations and therefore, the algorithm is more amenable to parallelization and distributed computing setting. This has also been noted by~\cite{chaillat2012faims}.


 To summarize the asymptotic cost of factorization using RandSVD at the leaf level is $\bigO(r_LMN + r_L^2(N +M))$ and using ppACA it is $\bigO( r_L^2(N + M))$. Numerical evidence suggests that there is compression at every level (and therefore $r_L \ll R$) and this justifies the use of this hierarchical approach.   

\subsection{Accuracy}\label{sec:acc} We now discuss the accuracy of the recursive low-rank approximation. In the algorithms described above, there are two sources of error - due to the low-rank truncation at the leaf level and the error accumulated due to the agglomeration process at all other levels in the tree. In order to analyze the accuracy of the recursive SVD computation, we first consider the accuracy of the agglomeration step at one level. We consider the matrix $\bH$ which has a partitioned as $\bH = [\bH_1^T, \bH_2^T]^T $. Suppose we compute a low-rank approximation to $\bH_i \approx \hat{\bH}_i\define \bU_i\bV_i^T$ for $i=1,2$ using the techniques described in Section~\ref{sec:compression}. We assume that the low-rank matrices satisfy the bounds $\| \bH_i -  \hat{\bH}_i\| \leq \varepsilon \| \bH_i\| $ for $i=1,2$. We can then bound the error in the approximation $\bH \approx \hat{\bH}\define \bU\bV^T$.

\begin{align*}
\norm{ \bH - \hat{\bH} }  \quad \leq & \quad  \norm{ \hat{\bH} - [\hat{\bH}_1^T,\hat{\bH}_2^T]} + \norm{\bH_1-\hat{\bH}_1} + \norm{\bH_2 - \hat{\bH}_2} \\ 
\leq & \quad \varepsilon \left(\norm{  [\hat{\bH}_1^T,\hat{\bH}_2^T]^T} \right)   + \varepsilon\left(\norm{\bH_1}+\norm{\bH_2}\right) \\  
\leq & \quad \varepsilon \left(\norm{  \hat{\bH}_1} +  \norm{\hat{\bH}_2} \right)   + \varepsilon\left(\norm{\bH_1}+\norm{\bH_2}\right) \\  
\leq & \quad (2\varepsilon + \varepsilon^2) \left(\norm{\bH_1}+\norm{\bH_2}\right)  = 2\varepsilon\left(\norm{\bH_1}+\norm{\bH_2}\right) + \bigO(\varepsilon^2) 
\end{align*}•
We have used the fact that the strategy that is used to truncate the rank for the agglomeration is the same as one to compute the low-rank compression. Furthermore, we have also have used the inequality that 
\[ \norm{  \hat{\bH}_i} = \norm{  \hat{\bH}_i - \bH_i + \bH_i } \leq \norm{  \hat{\bH}_i - \bH_i} + \norm{   \bH_i} \leq (1+ \varepsilon )\norm{\bH_i} \]
We now extend it to the case where $N_s > 2$ by recursively applying the error bound that was derived above 
\begin{align*}
\norm{ \bH - \hat{\bH} } \quad \leq & \quad \sum_{b \in \mathcal{L}(T_{I})} \norm{\bH_b - \hat{\bH}_b}  + \sum_{\ell=0}^{L-1}\sum_{b \in T^{(l)}_{I}  }  \norm{\hat{\bH}_b -\text{Agg}\{ \hat{\bH}_{b'}  : b' \in S(b))\}  } \\
\leq & \quad (L+1) \varepsilon \left( \sum_{b \in \mathcal{L}(T_{I})} \norm{\bH_b} \right)+ \bigO(\varepsilon^2) 
\end{align*}•
If the computations  were performed in the Frobenius norm, then using the Cauchy-Schwarz inequality we can conclude that  $ \norm{\bH - \hat{\bH} }_{F}  \lesssim 2^{L/2}(L+1)\varepsilon \norm{\bH}_{F}$. In order to derive an equivalent relationship for the $2$-norm, we use the following inequality
\[ \normtwo{\bH - \hat{\bH} }  \leq \norm{\bH - \hat{\bH} }_F \lesssim  2^{L/2}(L+1)\varepsilon \norm{\bH}_{F} \leq \sqrt{N_r} 2^{L/2}(L+1)\varepsilon \normtwo{\bH}\]
where $N_r = \min\{M,N\} $. Therefore, in order to achieve a desired relative tolerance $\varepsilon_d $, the tolerance that is used in the low-rank approximation and the agglomeration can be computed as $\varepsilon \sim \varepsilon_d/2^{L/2}(L+1)\sqrt{N_r}$ for the Recursive SVD. 

This error bound although locally optimal can result in a low-rank factorization that may be suboptimal in terms of compression. For this reason, we propose an additional step 
for compressing the low-rank factors $\hat\bH = \bU\bV^T$. 

\begin{enumerate}
\item Compute thin QR factorizations $\bQ_\bU\bR_\bU = \bU$ and $\bQ_\bV\bR_\bV = \bV$
\item Compute SVD $\bR_\bU\bR_\bV^T = \bU'\boldsymbol\Sigma'(\bV')^T$
\item Truncate $R$ singular values and return $\bU = \bQ_\bU\bU_r'$ and $\bV = \bQ_\bV\bV_r'\boldsymbol\Sigma_r'$
\end{enumerate}
This additional cost is $\bigO(R^2(M+N) + R^3)$ and may be beneficial when $R \ll \min\{M,N\}$.

The take away is that the cost of storage and matvecs with $\hat\bH$ is $\bigO(R(M+N))$, which is 
critical when we need to access and multiply with the estimate repeatedly in the course of the optimization for the image parameters.  We now describe that optimization problem.

\section{Reconstruction algorithms}\label{sec:recon}

The recovery of the shape of the tumor and the chromophore concentrations from diffuse optics measurements is an ill-posed inverse problem. The inverse problem can be stated as follows: Given a set of measurements $\by$ that measures the scattered field $\bphi_s$ at multiple detector locations $\br_d$ due to incident field $\bphi_i$ from multiple source locations and illuminated at several different wavelengths, recover the spatially varying perturbation of absorption $\Delta \mu_a(\br,\lambda)$ and the concentration of the chromophore species. Standard approaches to deal with ill-posedness introduce some kind of regularization, such as Tikhonov regularization. Here, we consider the parametric level setup approach proposed in~\cite{aghasi2011parametric} (abbreviated as PaLS) and subsequently applied to the diffuse optical tomography problem in~\cite{larusson2012parametric}.

We briefly review the PaLS approach for parameterizing the shape perturbation. The characteristic function $\chi(\br)$ defined in equation~\eqref{eqn:perturbation} is taken as the $\tau$-level set of a Lipschitz continuous function $\varphi (\br): \mathcal{D} \rightarrow \mathbb{R}$. Using $\varphi(\br)$, the characteristic function $\chi(\br)$ can be expressed as 
\begin{equation} \chi(\br) = H\left(\varphi(\br) -\tau \right) \qquad \varphi(\br) = \sum_{k=1}^{n_p} \alpha_k\psi\left(\beta_k\normr{\br-\bchi_k}\right) \label{eqn:pals}
\end{equation}
where $H(\cdot)$ is the Heaviside function. In practice, we use smooth approximations $H_\varepsilon$ of the Heaviside function $H$, and its derivative denoted by $\delta_\varepsilon$. We represent the function $\varphi(\br)$ parametrically as weighted combinations of basis functions $\psi(\cdot)$ and we have $\normr{\br} = \sqrt{\normtwo{\br}^2 + \nu^2}$ and $\nu > 0$ is a small parameter to ensure that $\varphi$ is differentiable. Several choices are available for $\psi$ such as polynomials and radial basis functions. Here we choose the compactly supported radial basis functions that were previously used in~\cite{aghasi2011parametric}. The coefficients $\alpha_k$ control the magnitude of the radial basis functions, $\beta_k$ control the width and $\bchi_k$ control the centers. The basis functions and their number control how fine or coarse the representation will be. On the one hand, having a large number of basis functions will be beneficial in reconstructing fine scale features, however, it has additional associated computational cost and further exacerbates the non-convexity.

The parameters that need to be estimated are collected in a vector $\bp = [\balpha^T,\bbeta^T,\bchi_x^T,\bchi_y^T,\bchi_z^T]^T \in \mathbb{R}^{N_p}$, where $N_p = (d+2)n_p$ and $d$ is the dimension of the problem. The reconstruction problem can now be stated as the minimization of the following functional
\begin{equation}\label{eqn:minimization}
  \hat{\bc}, \hat{\bp} \quad \define \quad \argmin_{\bc, \bp}  \text{ } \normtwo{\beps}^2 \text{ }=  \text{ }\normtwo{\bW(\by - \bD(\bp) \bc)}^2   
\end{equation}
where the columns of $\bD(\bp)$ are given by $\bE_i\bH\bmu(\bp)$ and $\bc = [c_1,\dots,c_{\nsp}]^T$ represent the concentration of the chromophores.  

\begin{algorithm}[!ht]
\begin{algorithmic}[1]
\STATE Given tolerances $\tau_1$ and $\tau_2$ and initial guess for PaLS parameters $\bp$
\WHILE {$\normtwo{\beps} \leq \tau_1$}
\STATE $\bc  =  \left(\bW\bD(\bp) \right)^{+} (\bW \by)$
\WHILE {$\normtwo{\beps} \leq \tau_2$}
\STATE\label{step:lm} $(\bJ^T\bJ + \nu \bI) \delta \bp = - \bJ^T \beps $ 
\\ \COMMENT{//The parameter $\nu$ is chosen by a Levenberg-Marquardt procedure. }
\STATE $\bp \leftarrow \bp + \delta \bp$
\ENDWHILE
\ENDWHILE
\RETURN Shape parameters $\bp$ and chromophore concentration $\bc$
\end{algorithmic}
\caption{Optimization procedure for solving shape parameters $\bp$ and chromophore concentrations $\bc$}
\label{alg:opt}
\end{algorithm}

The resulting optimization problem is solved by alternating between solving for the concentration parameters $\bc$ which is a linear least-squares problem and solving for the PaLS parameters $\bp$ using a Levenberg-Marquardt procedure. The optimization algorithm requires constructing the Jacobian 
\[ \bJ = \frac{\partial\beps}{\partial\bp} = - \bW\bar{\bE}\bH\frac{\partial\bmu}{\partial\bp} \]
where $\bar{\bE} = \sum_{i=1}^{\nsp}c_i\bE_i$. Analytical expressions for the derivatives $\frac{\partial\bmu}{\partial\bp}$ are provided in~\cite{aghasi2011parametric}. The stopping criteria for the iterative procedure is chosen according to the discrepancy principle, i.e., the iterations are terminated when the norm of the residuals is less than the noise level up to a user defined constant $\gamma > 1$. In mathematical terms, the stopping criterion  becomes $\normtwo{\beps} \leq \gamma \normtwo{\bEta}$ and $\bEta$ is the noise defined in Equation~\eqref{eqn:measurement}. More efficient algorithms are available for the reconstruction of PaLS parameters, for example, see TREGS~\cite{de2011regularized}. However, we have chosen the Levenberg-Marquardt algorithm for its relative simplicity of implementation.

Recall in Section~\ref{sec:compress}, we used a compressed low-rank representation of the measurement operator $\bH \approx \bU \bV^T$. Let us denote $\hat{\bH} = \bU\bV^T$ and the error as $\bE_\bH$ in the low-rank truncation process such that $\bH = \hat\bH +\bE_\bH $. From the results in Section~\ref{sec:compress} we know that $\norm{\bE_\bH} \leq \varepsilon \|\bW\bar\bE\|\norm{\bH} $ and define $\bar\varepsilon \define \varepsilon \|\bW\bar\bE\|$. To simplify the theoretical analysis, we rescale $\by \leftarrow \bW$, $\bH \leftarrow \bW\bar\bE\bH$ and $\hat{\bH}  \leftarrow \bW\bar\bE\hat\bH$.  The approximate Jacobian $\bar{\bJ}$ is now given by the expression $\bar{\bJ} = - \bar{\bH} \frac{\partial\bmu}{\partial\bp}$, so that we have $\bJ = \bar{\bJ} + \bE_\bJ $, where $\bE_\bJ =\bE_\bH\frac{\partial\bmu}{\partial\bp}$. It can be readily shown that the approximation to the objective function $\bbf(\bp) \define \normtwo{\by-\bH\bmu(\bp)}^2$ and the gradient $\nabla_\bp \bbf (\bp) \define -\frac{\partial\bmu}{\partial\bp}^T\bH^T(\by-\bH\bmu(\bp))$ (and the equivalent quantities $\hat\bbf$ and $\nabla_\bp \hat\bbf$ with the approximation $\bar\bH$ instead of $\bH$) satisfy the following approximation bounds
\begin{align*}
|\bbf-\hat{\bbf}|  \quad \leq & \quad  2\bar\varepsilon\normtwo{\by-\hat\bH\bmu} \normtwo{\bH}\normtwo{\bmu}+ \bigO(\bar\varepsilon^2) \\ 
\norm{\nabla_\bp \bbf- \nabla_\bp \hat{\bbf}}\quad \leq & \quad \bar\varepsilon \normtwo{\frac{\partial\bmu}{\partial\bp}} \normtwo{\bH} \left( \normtwo{\bmu}\normtwo{\bH} + \normtwo{\by-\hat\bH\bmu} \right) 
\end{align*}

Furthermore, assume that $\norm{\nabla_\bp \bbf- \nabla_\bp \hat{\bbf}} \leq \tau_g \norm{ \nabla_\bp \hat{\bbf}}$, then the acute angle $\theta$ between the gradient $\nabla_\bp \bbf$ and $\nabla_\bp \hat\bbf$ satisfies the following inequality
\[ \cos\theta \geq \frac{1-\tau_g^2}{\sqrt{1 + \tau_g^2}}\]

The result follows from the result in~\cite[Lemma 3.1]{yue2013accelerating}. We assume that the angle $\theta$ is acute, i.e., $\tau_g < 1$ which is always possible since we can control the error $\bar\varepsilon$ and therefore the tolerance $\tau_g$. We now present a result that bounds the error between the true and the perturbed search directions in step~\ref{step:lm} of Algorithm~\ref{alg:opt}.

\begin{propos}
Assume that the derivative $\frac{\partial\bmu}{\partial\bp}$ is full rank and $N_p \leq R$, where  $R$ is the effective rank of the low-rank representation $\hat\bH$ and $N_p$ is the number of PaLS parameters. In Algorithm~\ref{alg:opt} let $\delta \bp$ be the search direction corresponding to the exact Jacobian $\bJ$ and let $\delta\bar{\bp}$ be the search direction corresponding to the approximate Jacobian $\bar{\bJ}$. Then, we can bound the error between the two search directions as 
\begin{equation}
\normtwo{\delta \bp - \delta\bar{\bp}} \leq \left[ \eta(\nu) \normtwo{\delta\bar{\bp}} + \frac{\normtwo{\beps}}{\nu + \bar{\sigma}_{N_p}^2}\right] \normtwo{\bE_\bJ}
\end{equation}
where  the factor $\eta(\nu)= \max_{\bar{\sigma}_{N_p} \leq \sigma \leq \bar{\sigma}_1} \sigma/(\nu + \sigma^2)$ and $\bar\sigma_1$ and $\bar{\sigma}_{N_p}$ are upper and lower bounds for the singular values of the unperturbed Jacobian $\bJ$.

\end{propos}
The proof is readily obtained by an application of the result~\cite[Theorem 3.1]{ipsen2011rank}. The matrix $\frac{\partial\bmu}{\partial\bp}$ is full-rank and the dimension of $\bp$ denoted by $N_p$ is smaller than the rank of $\hat\bH$, therefore the Jacobian $\bJ$ and perturbed Jacobian $\bar{\bJ}$ are full-rank and satisfy the requirements of the theorem. If $R < N_p$ then the Jacobian is rank-deficient and we could consider a subset selection procedure similar to~\cite{ipsen2011rank}.

\section{Numerical Experiments}~\label{sec:results}
We present some results of the algorithms that we described in Section~\ref{sec:krylov}. For the rest of this section, we consider the following test problem. The geometry under consideration is a breast shaped phantom that is compressed between two flat plates (see Figure~\ref{fig:phantom}). At its widest it is 12 cm long and the maximum thickness is 5 cm. The domain is discretized using gmsh~\cite{geuzaine2009gmsh}, an open source 3D finite element mesh generator. The finite element matrices corresponding to the discretized representations of the partial differential equations given by equation~\eqref{eqn:matrices} are computed using FEniCS~\cite{logg2010dolfin} accessed using its Python interface. The boundary $\partial \Omega_R$ for which refractive index mismatch conditions are applied are assumed to be the flat top portions of the boundary, where as zero Dirichlet boundaries are applied on the rest of the boundary.

We consider the background medium to be composed of $\nsp = 4$ species, oxygenated and de-oxygenated hemoglobin denoted as HbO$_2$ and HbR respectively, water H$_2$O and lipids. These species have been specifically chosen since they are the most optically active chromophores, found in breast tissue, in the wavelength range $[600,1000]$ nm. The concentration of the various species in the background have been summarized in Table~\ref{tab:conc}. The extinction coefficients for the species have been found in the literature~\cite{prahl}. We take the value~\cite{grosenick2005timeb} of $\Psi = 9.4$ based on the wavelength of $600$ nm and the prefactor $b=1.4$ in Equation~\eqref{eqn:difflambda}.

\begin{table}[!ht]\centering
\begin{tabular}{|c|c|c|c|c|} \hline
Species & HbO$_2$ &HbR& H$_2$O & Lipids \\ \hline
Units & $\mu$M & $\mu$M & $\%$ & $\%$ \\ \hline
Background & 17 & 7  & 0.15 & 0.6 \\ \hline
Tumor & 25 &  15  & 0.25  & 0.5 \\ \hline
\end{tabular}
\caption{Concentration of different species in the background and the tumor. }
\label{tab:conc}
\end{table}

\subsection{Forward solver}	

For the preconditioner we choose an incomplete LU factorization implemented using SuperLU~\cite{lishao10} and we considered the parameters \verb;fill_factor; $\in \{5,10,15\}$ and \verb;drop_tol; $\in\{10^{-3},10^{-4},10^{-5}\}$. We consider the following transformation of the linear systems~\eqref{eqn:multipleshifted} 
\[ \left( \bK + \bar{\sigma}'\bR + \sigma_j\bM +( \sigma_j' -\bar{\sigma}')\bR  \right)\bx_j = \bb\]
for $j = 1,\dots,N_\lambda$ and $\bar{\sigma}'$ is the mean of $\sigma_j'$. We then define $\bK \leftarrow \bK + \bar{\sigma}'\bR $ and $\sigma_j' \leftarrow \sigma_j' - \bar{\sigma}'$. This transformation essentially the leaves the solution unchanged but improves the convergence of our solver, since the modified matrix $ \bK + \bar{\sigma}'\bR$ contains average information about the refractive index mismatch boundary conditions. Other transformations involving the minimum or maximum over $\sigma_j'$ may also be considered.  

We now report the results of our solver on a variety of test problem sizes, and preconditioner parameter such as fill factor and drop tolerance. The column labeled `Iter' reports the total number of iterations across $100$ wavelengths, `MVP [s]' reports the CPU time spent on matrix-vector products as well as application of the preconditioner, and finally `Tot. [s]' reports the total CPU run time of Algorithm~\ref{alg:auggmres} (including the pre-computation time for computing the initial guess $\tilde{\bx}_{0,j}$ and the matrices $\bU, \bC$ obtained by solving the shift-invariant system and generating the augmented space for each system $\bU_j$ and $\bC_j$). From the table it can be seen that by increasing the dimension of the deflation space $k$, the total number of iterations decrease but the cost per iteration increases as a result of extra orthogonalization w.r.t. $\bC_j$. Therefore, there is a trade-off between the number of iterations and total run time and adding additional vectors in the deflation space is a case of diminishing returns. For the range of parameters we experimented with, typically $k=5,10$ produces the best results in terms of total CPU time. However, for problem sizes larger than we are considering, the cost of matrix-vector products may be the dominant cost so that it might be beneficial to use a larger deflation space.

\begin{table}\centering
\begin{tabular}{|c|c|c|c|c|c|c|c|c|c|} \hline
\multicolumn{10}{|c|}{Varying problem size, drop tol $=10^{-4}$, fill  factor $=10$}\\\hline
\multirow{2}{*}{Dim k} & \multicolumn{3}{|c|}{$N=16,271$} & \multicolumn{3}{|c|}{$N=52,425$} & \multicolumn{3}{|c|}{$N=87,431$}\\ \cline{2-10}
•	& Iter 	& MVP [s] 	& Tot. [s] 	& Iter 	& MVP [s] & Tot. [s] 	& Iter & MVP [s] & Tot. [s] 
\\ \hline
$0$ 	& $2927$  & $27.90$ & $41.84$ 	& $3992$	& $146.39$ 	& $238.11$	& $5426$  & $388.35$	& $740.34$  \\ \hline
$5$ 	&  $2164$	& $18.09$	& $32.84$	& $3028$	& $124.39$	& $232.86$	& $4922$ 	& $315.79$ & $495.68$\\  \hline
$10$ & $2093$ 	& $18.94$	& $34.29$ 	& $2845$	& $116.44$	& $195.03$	& $4768$	& $300.87$& 	$469.59$	\\ \hline
$15$ & $2072$ 	& $19.34$	& $35.21$ 	& $2850$	& $103.59$	& $180.86$	& $4709$	& $311.00$ &  $479.408$	\\ \hline
\multicolumn{10}{|c|}{Varying fill factor, $N=52,425$, drop tol; $= 10^{-4}$  } \\ \hline
& \multicolumn{3}{|c|}{fill factor $= 5$} & \multicolumn{3}{|c|}{fill factor $= 10$} & \multicolumn{3}{|c|}{fill factor $= 15$}\\ \hline
$0$ 	& $5080$  & $129.91$ & $297.57$ 	& $3992$	& $146.39$ 	& $238.11$	& $3688$  & $195.37$  &  $262.99$  \\ \hline
$5$ 	&  $3710$	& $105.21$ & $222.64$	& $3028$	& $124.39$	& $232.86$	& $2809$ 	& $165.09$  & $248.91$\\  \hline
$10$ & $3614$ 	& $98.20$	& $210.50$ 	& $2845$	& $116.44$	& $195.03$	& $2664$	& $145.67$  & 	$218.17$	\\ \hline
$15$ & $3729$ 	& $91.93$	& $197.10$ 	& $2850$	& $103.59$	& $180.86$	& $2651$	& $131.70$  &  $194.64$	\\ \hline
\multicolumn{10}{|c|}{Varying drop tol, fill factor = $10$, $N=52,425$ } \\ \hline
& \multicolumn{3}{|c|}{drop tol $= 10^{-3}$} & \multicolumn{3}{|c|}{drop tol $= 10^{-4}$} & \multicolumn{3}{|c|}{drop tol $= 10^{-5}$}\\ \hline
$0$ 	& $3485$  & $99.61$ & $179.86$ 	& $3992$	& $146.39$ 	& $238.11$	& $5293$  & $193.75$    & $314.37$  \\ \hline
$5$ 	&  $2665$	& $77.21$	& $169.98$	& $3028$	& $124.39$	& $232.86$	& $3859$ 	& $143.45$  &   $248.03$\\  \hline
$10$ & $2541$ 	& $75.81$	& $161.72$ 	& $2845$	& $116.44$	& $195.03$	& $3721$	& $155.99$  & 	$259.73$	\\ \hline
$15$ & $2505$ 	& $71.76$	& $148.10$ 	& $2850$	& $103.59$	& $180.86$	& $3744$	& $144.06$ &  $228.62$	\\ \hline
\end{tabular}•
\caption{Summary of augmented GMRES solver for different problem sizes, and preconditioner parameters. The number $k$ refers to the dimension of the deflation space. All systems were solved till it converged to a relative tolerance of $10^{-8}$. As can be seen, on average $10-20\%$ improvement was observed for all systems  in terms of total computation time, by using deflation. }
\label{tab:timingforward}
\end{table}

We would like to emphasize that care should be taken to interpret the results in Table~\ref{tab:timingforward}. The algorithm has been implemented in Python (which is an interpreted language) and because it uses pre-compiled code for parts of the computation, the timing results may be slightly different if the entire algorithm were implemented in a single programming language. In particular, we expect the overall computation time would be lower if using a compiled language such as C/C++ and the performance gains from our algorithm to be higher.

\subsection{Compression}
We now discuss the results of the compression scheme presented in Section~\ref{sec:compress}. We consider the same geometry that was used in the previous subsection. The sources are placed on the top of the phantom whereas the receivers are placed on the bottom. For each source, there are $9$ detectors constrained to move along with the source that are evenly placed co-axially with the source with $0.5$ cm distance from each other. The number of sources varied from $N_s = 4,\dots,25$ and the number of wavelengths vary between $N_\lambda = 11,\dots,81$. With these parameters, the maximum number of parameters are $18,225$. The discretized grid has  $N = 52,425$ degrees of freedom.

\begin{figure}[!ht]
\centering
\includegraphics[scale=0.33]{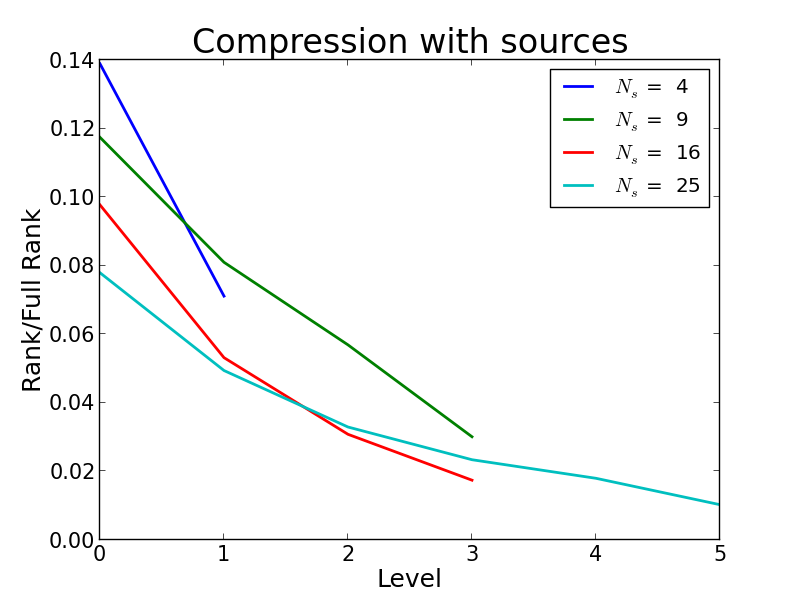}
\includegraphics[scale=0.33]{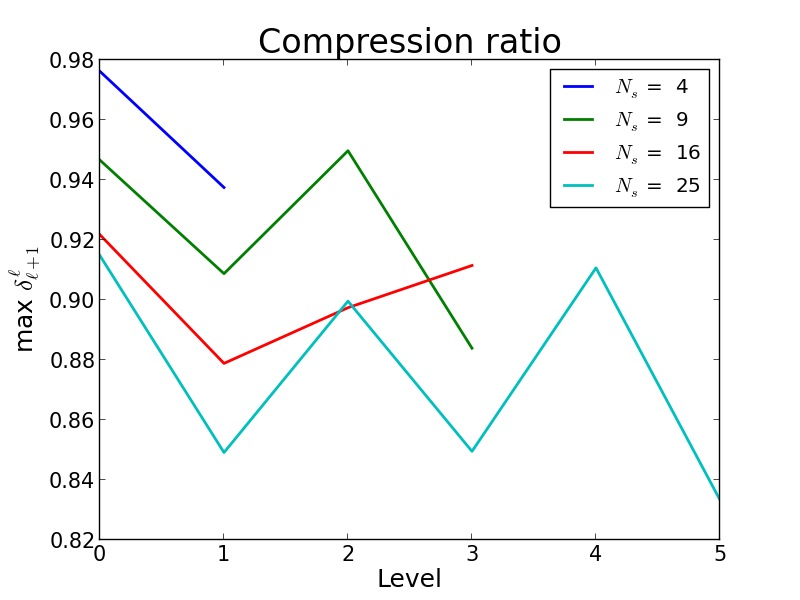}
\caption{(left) Compression defined as Rank/Full rank of the recursive SVD computed as a function of the level of the tree, with increasing number of sources. Here Rank is defined as the number of columns of $\bU$ and Full Rank $=\min\{M,N\} $. (right) The compression ratio $\delta^{\ell}_{\ell+1} \define r_\ell/(2r_{\ell+1})$ as function of the level. As can be seen, on average, there is higher levels of compression at higher levels in the tree. Here, level $0$ is the root of the tree. A tolerance of $10^{-6}$ was used for the truncation. Furthermore, $N_{ds} =9$ and $N_\lambda = 41$.}
\label{fig:source}
\end{figure}
In the examples, we will describe, we have used the randomized SVD for both computing the low-rank factorization at the leaf level and to compute the agglomeration of the low-rank factors as we go up the tree. Similar results are obtained using the partially pivoted Adaptive Cross Approximation and will not be displayed here. We first consider the compression by varying the number of sources and keeping all other parameters fixed. We assume that $N_{ds} =9$ and $N_\lambda = 41$ and a tolerance of $10^{-6}$ was used for truncating the rank of the sub-blocks and the agglomeration. As can be seen, at higher levels in the tree (closer to the root) we observe a higher level of compression because there is a greater redundancy of information globally as opposed to locally. Furthermore, with increasing number of sources we observe a higher level of compression which implies that there is redundancy both in terms of wavelengths and the source-detector positions. The results are displayed in Figure~\ref{fig:source}. We also plot the maximum compression ratio $\delta^{\ell}_{\ell+1} \define r_\ell/(2r_{\ell+1})$ (computed across all the nodes at level $\ell$) which is the ratio of the ranks at level $\ell$ obtained by agglomerating 2 sub-blocks at level $\ell +1$ with ranks $\ell+1$. As can be seen,  $\delta^{\ell}_{\ell+1} < 1$ at all levels indicating that there is compression, not only at the leaves, but compression at every level in the tree. This justifies using a hierarchical compression scheme and the cost analysis performed in Section~\ref{sec:costs}.

\begin{figure}[!ht]
\centering
\includegraphics[scale=0.33]{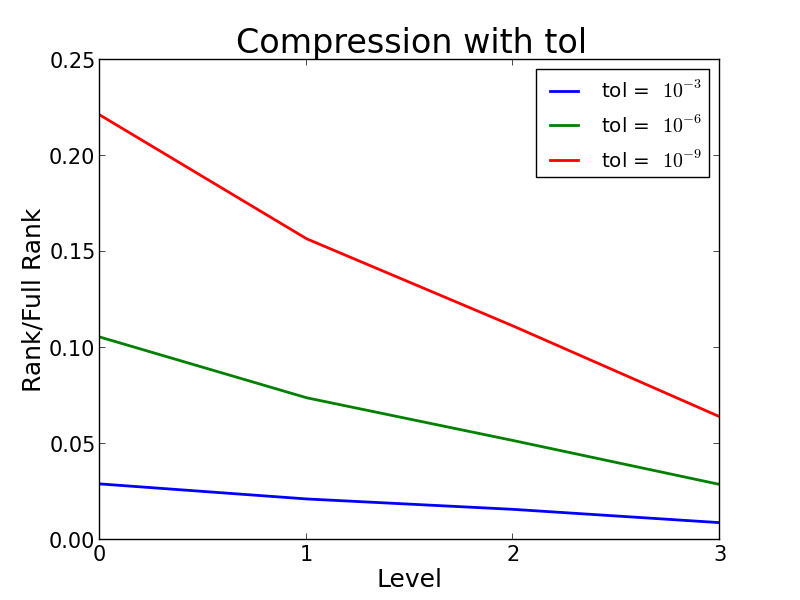}
\includegraphics[scale=0.33]{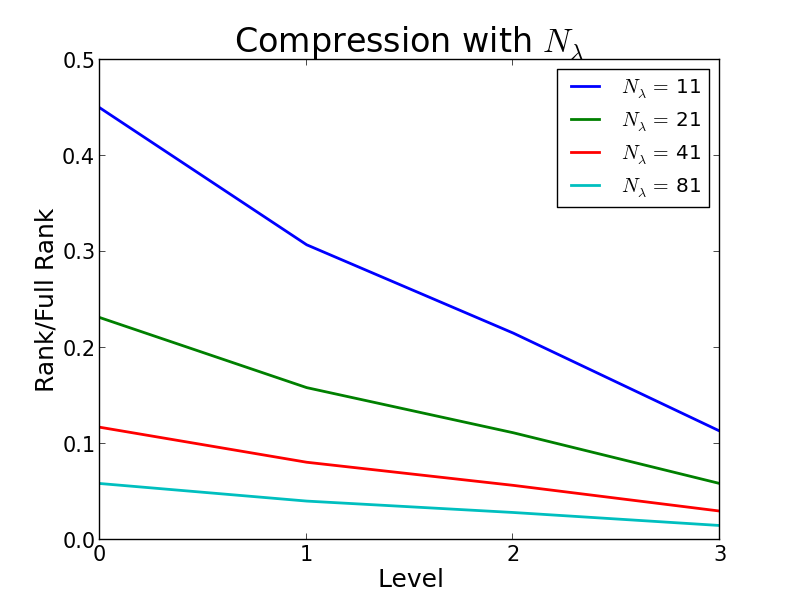}
\caption{(left) Compression defined as Rank/Full rank of the recursive SVD computed as a function of the level of the tree, with varying tolerance used to compress the low-rank factors. Here $N_{ds}=9$, $N_\lambda = 41$ and $N_s = 9$.   (right) The compression of the recursive SVD computed as a function of the level of the tree, with varying $N_\lambda$. Here $N_{ds}=9$, tol = $10^{-6}$ and $N_s = 9$. }
\label{fig:tollambda}
\end{figure}

Next we compute the compression as a function of tree level with varying tolerance used to truncate the ranks of the sub-blocks. The results are presented in Figure~\ref{fig:tollambda}. As can be seen with a higher tolerance the ranks at each level decreases dramatically. Here, we have fixed  $N_{ds}=9$, $N_\lambda = 41$ and $N_s = 9$. We also present results of computing the compression as a function of tree level with varying number of wavelengths used to illuminate the object. All other parameters are fixed as $N_{ds}=9$, tol = $10^{-6}$ and $N_s = 9$. The results are also presented in Figure~\ref{fig:tollambda}. We can see that with increasing number of wavelengths there is a higher level of compression at higher levels in the tree.

 Finally, we compare the run time of the recursive SVD algorithm proposed in Section~\ref{sec:compress} with RandSVD applied to the entire measurement operator. As can be seen, the computational time for the recursive algorithm is far lower and is therefore, more efficient. 

\begin{figure}[!ht]
\centering
\includegraphics[scale=0.3]{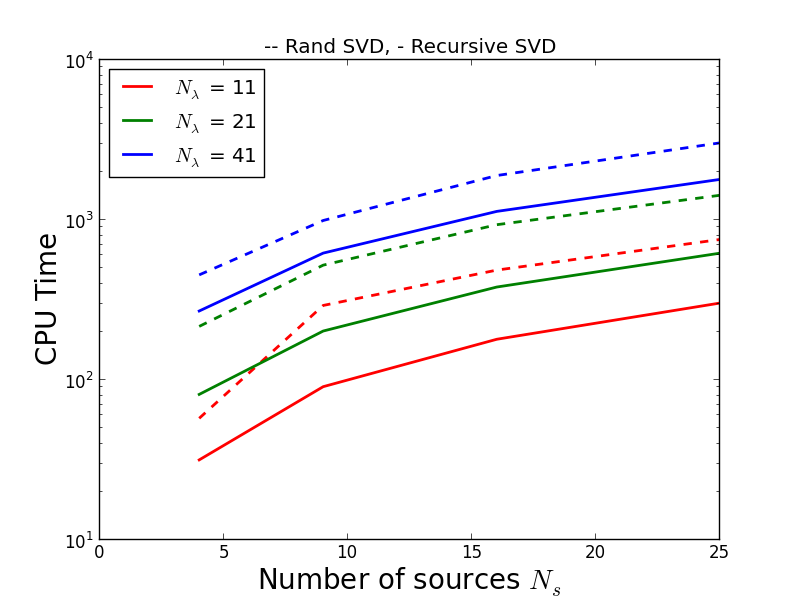}
\caption{Comparison of time taken to factorize the Born operator $\bH$ to a tolerance of $10^{-6}$ with $N_s $ ranging from $4$ to $25$ and $N_{ds}$ fixed at 9. We compare the CPU timing of RandSVD in Algorithm~\ref{alg:randsvd} applied to the entire matrix directly with the Recursive SVD algorithm proposed in Section~\ref{sec:compress}. The time for constructing $\bH$ is not included. The grid size was $16,721$.}
\label{fig:comptime}
\end{figure}

\subsection{Reconstruction results}

In the following experiments we will use the following metrics to measure the error in the shape perturbation. In Experiment $1$ and $2$ we use the full Born operator to generate the measurements, however, we use the compressed operator in the reconstructions. In Experiment $3$ we use measurements from the full diffusion equation and use the compressed Born model in the reconstruction. As a result, we avoid Let $\bmu$ denote the true shape perturbation and let $\hat\bmu$ denote the recovered shape perturbation. The first is standard relative L$_2$ error $\normtwo{\bmu-\hat\bmu}/\normtwo{\bmu}$. For piecewise constant medium, a different metric known as the Dice coefficient (see for example,~\cite{larusson2012parametric}) provides more information about localization of the perturbation. It can be defined as 

\[ D(\bmu,\hat\bmu) \define 2 \frac{| \bmu \bigcap \hat\bmu|}{|\bmu| + |\hat\bmu|} \]
where $| \bmu \bigcap \hat\bmu|$ corresponds to the number of non-zero pixels that are contained in both the true shape $\bmu$ and the reconstructed shape $\hat\bmu$ whereas $|\bmu|$ and $|\hat\bmu|$ correspond to the number of non-zero pixels in the true and the reconstructed shapes respectively. In order to If the object is recovered fully, then the Dice coefficient would be $1$. 

\textbf{Experiment 1}: In this experiment we study how the accuracy of the low-rank representation for $\bH$ affects the error in the reconstruction. We consider $N_s=4$ sources and $N_{ds} = 9$ detectors/source and $N_\lambda=25$ wavelengths totaling $900$ measurements generated using the Born model but with the full measurement operator $\bH$. We only consider reconstruction of the shape perturbation and consider the concentration of the chromophores as known and provided in Table~\ref{tab:conc}. Gaussian noise is added such that the signal-to-noise (SNR) ratio is $33$ dB to simulate observational noise. The SNR is defined as
\[ \text{SNR} = 20\log_{10} \frac{\normtwo{\by}}{\normtwo{\boldsymbol\eta}}\]
 where the noise $\bEta$ has been defined in Equation~\eqref{eqn:measurement}. The grid size is chosen to be $52,425$. A ``true'' shape perturbation is obtained by using three randomly generated basis functions which leads to $N_p = 15$. We report the rank of the measurement operator as a function of the global tolerance used for low-rank representation and the error in the reconstruction using the low-rank operator $\hat\bH$. The results are reported in Table~\ref{tab:exp1}. As can be seen that for very small tolerance, the error in the reconstruction is hardly noticeable. Therefore, a larger rank of the low-rank representation $\hat\bH$ does not affect the reconstruction error below a certain tolerance and therefore the compressed operator $\hat\bH$ can be used as a surrogate for the full matrix $\bH$ with little or no loss in accuracy in the reconstruction. When the number of measurements is large, this can represent significant savings in computational time. 

\begin{table}[!ht]\centering
\begin{tabular}{|c|c|c|c|c|} \hline
tol  & $10^{-3}$ & $10^{-6}$ & $10^{-9}$\\ \hline
Rank & 48&  120& 288 \\ \hline
L$_2$err. &  $40.62 \%	$ & $40.52 \%$&  $40.52 \%$ \\\hline
Dice &$0.82$ & $0.827$ & $0.827$  \\ \hline
Time [s]& $2.20$ & $2.23$ & $2.32$ \\  \hline
\end{tabular}
\caption{Rank of the compressed operator and the error in the reconstruction of the shape perturbation as a function of tolerance used to compress $\bH$. See experiment 1 for more details. `Time [s]' indicates the CPU time of solving the optimization problem with the compressed operator. See Experiment 1 for more details.}
\label{tab:exp1}
\end{table}

\textbf{Experiment 2}: In this experiment we study the reconstruction of chromophore concentrations as well as the shape perturbation. The number of measurements and the process of generating them is the same as Experiment 1 except with SNR $30$ dB. For the reconstruction, the radial basis functions were randomly initialized and a truncated measurement operator $\hat\bH$ computed using tolerance $10^{-6}$ was used in the reconstruction. The error of the reconstruction of the chromophore concentrations as well the relative L$_2$ error of the shape perturbation are reported in Table~\ref{tab:exp2}. As can be seen from the Table, the added difficulty in recovering the shape perturbation as well as the chromophore concentrations affects the reconstruction error of the shape perturbation slightly. Moreover, the concentrations of the chromophore species are recovered fairly accurately. This is consistent with the observations in~\cite{larusson2013parametric,larusson2011hyperspectral,larusson2012parametric}.

\begin{table}[!ht]\centering
\begin{tabular}{|c|c|c|c|c|c|c|} \hline
Species & HbO$_2$ &HbR& H$_2$O & Lipids & L$_2$ err. & Dice \\ \hline
 Recon. & $1.6 \%$ & $0.2 \%$ & $4.9 \%$ & $4.4\%$ &  $52.72 \%$ & $0.79$ \\\hline
\end{tabular}
\caption{Error in the reconstruction of different species and the shape perturbation corresponding to Experiment 2. }
\label{tab:exp2}
\end{table}
\begin{figure}
\centering
\includegraphics[scale=0.182]{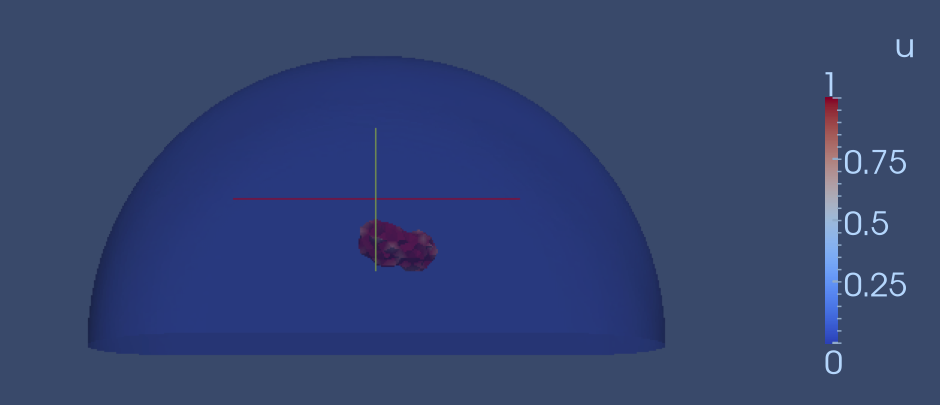}
\includegraphics[scale=0.177]{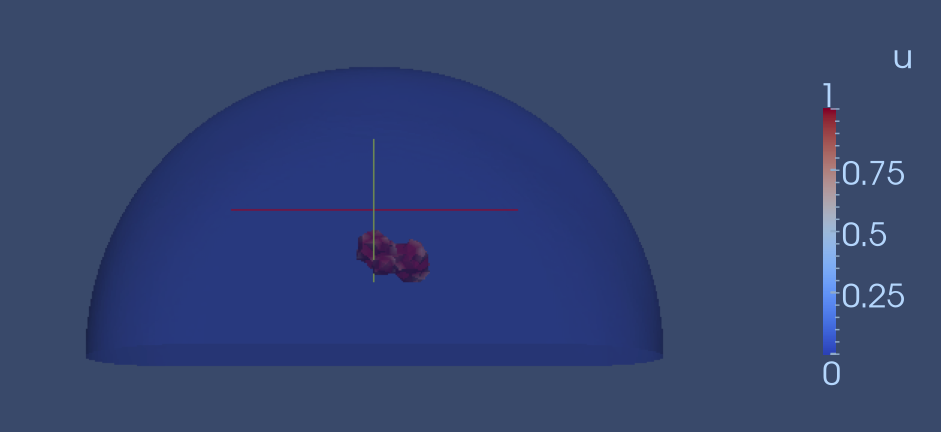}
\caption{Reconstruction in Experiment 2 (left) true anomaly and (right) reconstruction. The error is described in Table~\ref{tab:exp2}.}
\end{figure}

\begin{figure}[!ht]
\centering
\includegraphics[scale=0.3]{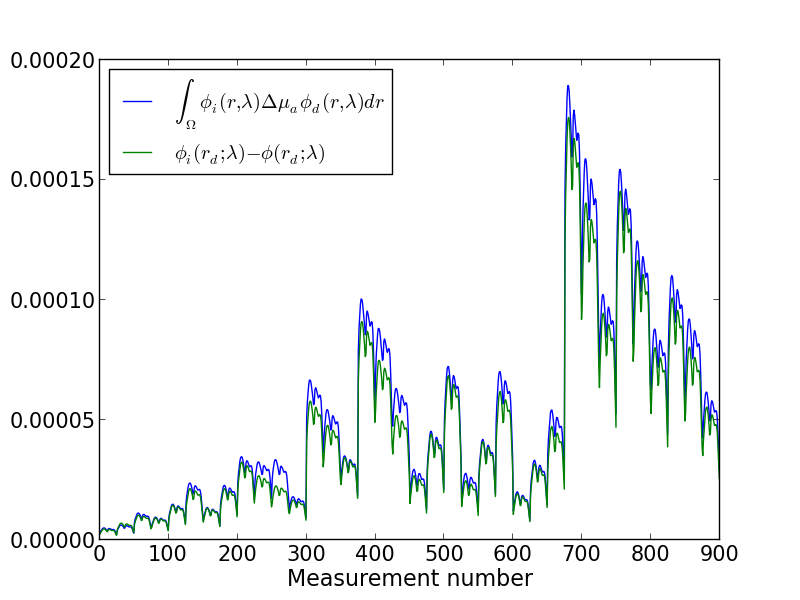}
\includegraphics[scale=0.22]{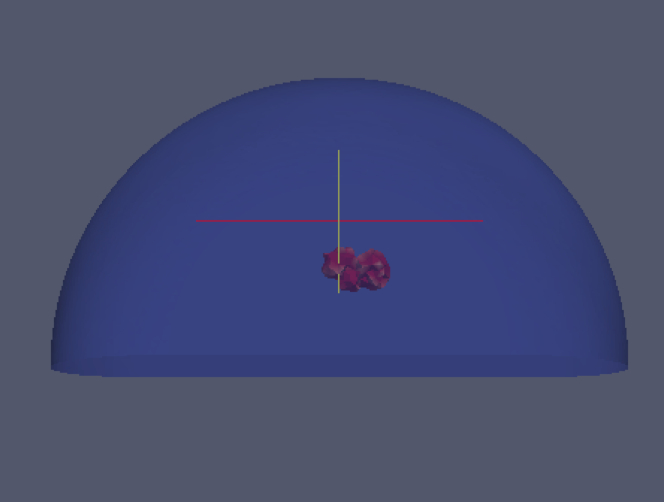}
\caption{(left) Comparison of the measurements generated using the Born model and that generated from the full diffusion equation. As can be seen, the agreement between the measurements is quite close. (right) Reconstruction of the shape perturbation using measurements from full diffusion equation. The errors are reported in Table~\ref{tab:exp3}.}
\label{fig:exp3}
\end{figure}

\textbf{Experiment 3}: In this example we examine the validity of the Born approximation. As mentioned earlier, experimental studies have validated the utility of the Born approximation (see for e.g.~\cite{larusson2012parametric}). However, in order to prevent committing an ``inverse crime'' we use data generated from the full diffusion equation, Equation~\eqref{eqn:diffusion} solving using finite elements with the same boundary conditions and use it to reconstruct both the chromophore concentrations and the shape parameter. The comparison between the measurements generated using the Born model and those generated from the full diffusion equation is provided in Figure~\ref{fig:exp3}. In addition a noise of $50$ dB was added to simulate observational error in realistic conditions. All other settings were the same as that in Experiment 2. The error of the reconstruction of the chromophore concentrations as well the L$_2$ error of the shape perturbation are reported in Table~\ref{tab:exp3}. The error in the reconstruction is higher than that obtained from Experiment 2. Since the full diffusion equation was used in generating the measurements, there is a modeling error which corresponds to about $20$ dB SNR and because of this, the Born model cannot exactly recover the shape perturbation. However, the reconstruction is still satisfactory as can be seen from the Figure~\ref{fig:exp3}.

\begin{table}[!ht]\centering
\begin{tabular}{|c|c|c|c|c|c|c|} \hline
Species & HbO$_2$ &HbR& H$_2$O & Lipids & L$_2$ err. & Dice \\ \hline
 Recon. & $5.1 \%$ & $ 1.2 \%$ & $1.2 \%$ & $9.3\%$ &  $71.80 \%$ & 0.65\\\hline
\end{tabular}
\caption{Error in the reconstruction of different species and the shape perturbation corresponding to Experiment 3. Measurements were generated from the full diffusion equation, i.e., Equation~\eqref{eqn:diffusion}. See also Figure~\ref{fig:exp3}.}
\label{tab:exp3}
\end{table}

\section{Conclusions and future work}
We have presented a fast algorithm for recovering shape of the perturbations and chromophore concentrations that is composed of three parts - a fast Krylov subspace approach for accelerating the solution of the incident and the adjoint field for multiple wavelengths, constructing a low-rank approximation to the sensitivity matrix $\bH$  using an approach that accounts for redundancies across wavelengths which is  then recursively combined across source-detectors pairs, and finally solving an optimization problem for recovering parameters with the low-rank approximation to $\bH$. The results indicate that significant gains can be obtained both in terms of computational costs and storage requirements.  We provide detailed numerical experiments that validates our claims and further provide a detailed analysis of the computational costs and the error. The algorithms were demonstrated on a challenging synthetic inversion case on a complex geometry which highlights the flexibility of our algorithms.

Future work includes extension of our algorithms to the fully nonlinear setting. In this setting, computation of the incident $\phi_i$ and adjoint fields $\phi_d$ and the construction of the measurement operator $\bH$ for the Born approximation, must be performed at every Newton or Gauss-Newton step. We therefore, believe that adopting the ideas proposed in this paper will be even more beneficial in the fully nonlinear case. Other possible extensions include a recycling strategy of the kind described in Section~\ref{sec:krylov} to multiple wavelengths and multiple right hand sides and to the full diffusion equation, Equation~\eqref{eqn:diffusion}. Additional work is currently underway in our lab to compare the reconstructions obtained the synthetic case with real data obtained from experiments. Based on previous work~\cite{larusson2013parametric,larusson2011hyperspectral,larusson2012parametric}, we anticipate that the reconstructions will indeed be excellent even in the hyperspectral case.

\section{Acknowledgements}
We would like to thank Nishanth Krishnamurthy, Pami G. Anderson, Jana Kainerstorfer, and Angelo Sassaroli for useful discussions. The first author would also like to thank Tania Bakhos with her help in generating the meshes. This work was supported by NIH Grant R01-CA154774. Additionally, the second author was supported by NSF Grant DMS 1217161.

\appendix

\section{Computing low-rank representations}\label{sec:lowrank}

\subsection{Randomized SVD}~\label{sec:randsvd}
Randomized algorithms for matrix decompositions were derived in a series of papers~\cite{liberty2007randomized,martinsson2011randomized,halko2011finding}.  In this work, we consider the algorithms described in~\cite{halko2011finding}, which was also used by~\cite{chaillat2012faims}. Suppose we wish to compute the rank $r$ decomposition of the matrix. The algorithm begins by computing a matrix $\bQ$ that approximates the column space of $\bA$ as 
\begin{equation}\label{eqn:rand} \| \bQ\bQ^*\bA  - \bA \| \leq \epsilon\end{equation}
where $\epsilon$ is a  user-defined tolerance. The matrix $\bQ$ is  obtained by computing the matrix-vector products of $\bA$ with the matrix $\boldsymbol\Omega_1 \in \mathbb{R}^{n\times(r+p)}$, with i.i.d. entries drawn from a standard normal distribution, $\mathcal{N}(0,1)$, and then computing a basis for the resulting matrix using QR or the SVD. Here, $p$ is an oversampling parameter that is chosen a priori and typically $p\sim 20$. A matrix $\bQ$ that satisfies the bound~\eqref{eqn:rand} can be converted into a low-rank representation using matrix manipulations. A discussion of the choice of oversampling factor and the low-rank conversion is described in~\cite{halko2011finding}.

The algorithm has good performance when the singular values of the matrix $\bA$ decay rapidly. If this is not the case, the power method is applied to improve the convergence of the algorithm. Since the rank of the matrix is not known a priori, we use an adaptive approach to estimate the range that is combined with an error estimator described in~\cite{halko2011finding}. The algorithm is summarized in~\ref{alg:randsvd}. The dominant computational cost is computing the matrix-vector products with the matrix $\bA$ which costs $\bigO(mn(r+p))$. The total cost of computing the low-rank representation is $\bigO\left(mn(r+p) + r^2(m+n)\right)$. 

\begin{algorithm}[!ht]
\begin{algorithmic}[1]
\STATE $r= 1$ (initial guess for rank), $k = 10$ (error estimator), $p = 20$ (oversampling factor)
\WHILE {\TRUE }
\STATE Contruct matrices $\boldsymbol\Omega_1 \in \mathbb{R}^{n\times (r+p)}$ and $\boldsymbol\Omega_2\in \mathbb{R}^{n\times k}$ with entries drawn from i.i.d. normal distribution $\mathcal{N}(0,1)$.
\STATE $\bY = \bA\boldsymbol\Omega_1$ 
\STATE $[\bU_r,\bS_r,\bV_r]$ = SVD$(\bY)$
\STATE $\bQ = \bU(:,1:r)$
\STATE $\bB\define\bQ^T \bA$
\STATE Error estimate  $e_r = \| \bA\boldsymbol\Omega_2 - \bQ\bQ^T \bA\boldsymbol\Omega_2\| $ 
\IF {$e_r > \varepsilon \bS_r(1,1)$ }
\STATE Increment $r$
\ELSE
\STATE  $[\bU_1,\bS,\bV]=$SVD($\bB $)
\STATE \textbf{Break}
\ENDIF
\ENDWHILE
\RETURN [$\bU , \bS, \bV] $ where $\bU=  \bU_r \bU_1 $.
\end{algorithmic}
\caption{Randomized SVD~\cite{halko2011finding}}
\label{alg:randsvd}
\end{algorithm}

\subsection{Adaptive Cross Approximation}
The idea behind the cross approximation is based on the result described in \cite{bebendorf2000approximation}, which states that supposing a matrix ${A}$ is well approximated by a low-rank matrix, by a clever choice of $k$ columns indexed as $\mathcal{J}$ and $k$ rows indexed as $\mathcal{I}$ of the matrix ${\bA}$, we can approximate $\hat{\bA}$ of the form

\[ \lVert {\bA} - \hat{{\bA}} \rVert \leq \varepsilon \qquad \hat{\bA} = {\bA(:,\mathcal{J})\bA(\mathcal{I},\mathcal{J})^{-1}\bA(\mathcal{I},:)} \]

 This decomposition relies on a result from~\cite{goreinov1997theory}, which states that if there is a sufficiently good low rank approximation to a matrix, then there exists a cross-approximation with almost the same approximation quality.  

\begin{algorithm}[!ht]
 \begin{algorithmic}[1]
  \STATE Initialize 
 \[ {\bR}_0 = {\bA},\qquad {\bS} = {\bzero}\]
  \FORALL {$ k = 0,1,2,\dots$} 	
  \STATE $(i_{k+1}^*,j_{k+1}^*) := \arg\max_{i,j}\lvert ({\bR}_k)_{ij} \rvert$ and $\gamma_{k+1}=\left({\bA}_{i^*_{k+1},j^*_{k+1}}\right)^{-1}$.
  \IF {$\gamma_{k+1} \neq 0$}
      \STATE Compute column ${\bu}_{k+1} := \gamma_{k+1}{\bR}_k{\be}_{j_{k+1}}$ and row ${\bv}_{k+1} := {\bR}_k^T{\be}_{i_{k+1}}$
      \STATE New residue and approximation
          \[ {\bR}_{k+1} := \bR_k -  \bu_{k+1}{\bv}^T_{k+1}\qquad {\bS}_{k+1} := \bS_k +  {\bu}_{k+1} {\bv}^T_{k+1}\]
  \ELSE \STATE Terminate algorithm with exact rank $k-1$
  \ENDIF
  \ENDFOR
 \end{algorithmic}
\caption{Cross Approximation using Full Pivoting~\cite{bebendorf2003adaptive}}
\label{alg:fullpivoting}
\end{algorithm}

Algorithm~\ref{alg:fullpivoting} describes a simple heuristic to compute such a cross approximation that is based on successive approximations by rank-$1$ matrices. It has the property that if the matrix ${\bA}\in\mathbb{R}^{m\times n}$ has an exact rank $r < \min\{m,n\}$, this algorithm will terminate in $r$ steps and defining 
\[
{\bS}_r= \sum_{k=1}^r{\bu}_k{\bv}_k^T 
\]
we have that $\bS_r={\bA}$ in exact arithmetic. Furthermore, it exactly reproduces the $r$ pivot rows and columns of $\bA$. Of course, the principal disadvantage of this algorithm is that, to generate a rank-$k$ approximation, it requires ${\cal{O}}(rmn)$ steps, which is not feasible for large matrices. The bottleneck arises from calculating the pivot indices $(i^*_k,j^*_k)$ which requires generating all the entries of the matrix ${\bA}$.

Several heuristic strategies have been proposed to reduce the complexity of the fully pivoting cross approximation algorithm. In particular, one such algorithm is called partially pivoted Adaptive Cross Approximation algorithm that  has a complexity $\bigO(r^2(m+n))$. A practical version of the algorithm, which includes a termination criteria based on an heuristic approximation to the relative approximation in the Frobenius norm, can be found in~\cite{bebendorf2003adaptive}. This is the version we will use in the rest of the paper.

\bibliographystyle{plain}
\bibliography{writeup} 
\end{document}